\documentclass[a4paper]
{cedram-smai-jcm}

\usepackage{tikz}
\usepackage{subcaption}

\usepackage[utf8]{inputenc}
\usepackage{array}
\usepackage{algorithm}
\usepackage{amsmath, amsthm, amssymb}
\usepackage[english]{babel}
\usepackage[inline]{enumitem}
\usepackage{booktabs}
\usepackage{mathtools}
\usepackage[numbers,sort&compress]{natbib}

\usepackage{tikz-cd}
\usetikzlibrary{shapes,arrows}
\usetikzlibrary{arrows.meta}
\usepackage{multicol}

\newcommand{\bs}{\mathbf{s}}
\newcommand{\Fs}{\mathcal{F}_*}

\newcommand{\bt}{\mathbf{t}}

\newcommand{\bn}{\mathbf{n}}
\newcommand{\bhn}{\hat{\mathbf{n}}}
\newcommand{\be}{\mathbf{e}}
\newcommand{\bhe}{\hat{\mathbf{e}}}
\newcommand{\bx}{\mathbf{x}}
\newcommand{\bhx}{\hat{\mathbf{x}}}

\newcommand{\sh}{\hat{s}}
\newcommand{\htau}{\hat{\tau}}
\newcommand{\br}{\mathbf{r}}
\DeclareMathOperator{\Div}{div}
\DeclareMathOperator{\curl}{curl}
\DeclareMathOperator{\rot}{rot}
\DeclareMathOperator{\dev}{dev}
\DeclareMathOperator{\airy}{airy}
\DeclareMathOperator{\alt}{Alt}
\DeclareMathOperator{\supp}{supp}

\newcommand{\jump}[1]{\ensuremath[\![#1]\!]}

\newcommand{\dx}{\ensuremath\,\text{d}x}
\newcommand{\dxhat}{\ensuremath\,\text{d}\hat{x}}
\newcommand{\ds}{\ensuremath\,\text{d}s}
\newcommand{\dshat}{\ensuremath\,\text{d}\hat{s}}
\newcommand{\hdiv}{H(\mathrm{div})}
\newcommand{\hcurl}{H(\mathrm{curl})}
\newcommand{\Fdiv}{\mathcal{F}^{\mathrm{div}}}
\newcommand{\Fdivdiv}{\mathcal{F}^{\mathrm{div}, \mathrm{div}}}
\newcommand{\tr}{\mathrm{tr}}
\newcommand{\vertiii}[1]{{\left\vert\kern-0.25ex\left\vert\kern-0.25ex\left\vert #1\right\vert\kern-0.25ex\right\vert\kern-0.25ex\right\vert}}
\newcommand{\verth}[1]{\vertiii{#1}_h}

\usepackage{todonotes}

\usepackage[normalem]{ulem}

\title[Piola-mapped elements]{Transformations for Piola-mapped elements}

\author[F.~R.~A.~Aznaran]{\firstname{Francis} \middlename{R.~A.} \lastname{Aznaran}}
\address{Mathematical Institute, University of Oxford, Oxford, OX2 6GG, UK}
\email{aznaran@maths.ox.ac.uk}
\thanks{The first author is supported by the EPSRC Centre for Doctoral Training in Partial Differential Equations: Analysis and Applications [grant number EP/L015811/1] and The MathWorks, Inc. The second author is supported by EPSRC grants EP/V001493/1 and EP/R029423/1. The third author is supported by National Science Foundation grant 1912653.}

\author[P.~E.~Farrell]{\firstname{Patrick} \middlename{E.} \lastname{Farrell}}
\address{Mathematical Institute, University of Oxford, Oxford, OX2 6GG, UK}
\email{patrick.farrell@maths.ox.ac.uk}

\author[R.~C.~Kirby]{\firstname{Robert} \middlename{C.} \lastname{Kirby}}
\address{Department of Mathematics, Baylor University; Sid Richardson Science Building; 1410 S.~4th St.; Waco, TX 76706.}
\email{robert\_kirby@baylor.edu}

\keywords{Finite element method, Piola transform, pullback, linear elasticity, Stokes, reference basis, Firedrake}
  
\subjclass{65N30, 65F08}

\begin{document}

\begin{abstract}
    The Arnold--Winther element successfully discretizes the Hellinger--Reissner variational formulation of linear elasticity; its development was one of the key early breakthroughs of the finite element exterior calculus. Despite its great utility, it is not available in standard finite element software, because its degrees of freedom are not preserved under the standard Piola push-forward. In this work we apply the novel transformation theory recently developed by Kirby~[\textit{SMAI-JCM}, 4:197--224, 2018] to devise the correct map for transforming the basis on a reference cell to a generic physical triangle. This enables the use of the Arnold--Winther elements, both conforming and nonconforming, in the widely-used Firedrake finite element software, composing with its advanced symbolic code generation and geometric multigrid functionality. Similar results also enable the correct transformation of the Mardal--Tai--Winther element for incompressible fluid flow. We present numerical results for both elements, verifying the correctness of our theory.
\end{abstract}

\maketitle
\vspace{-3mm}\section{Introduction}
\vspace{-2mm}

The use of a \emph{reference element} is critical in the evaluation of finite element basis functions and their derivatives.  One constructs, by some means, the basis functions on a particular fixed element and obtains the basis functions on an arbitrary cell $K$ through a mapping.  For scalar-valued function spaces, this is frequently a simple pullback (change of coordinates).  However, vector-valued elements discretizing $\hdiv$ or $\hcurl$, as well as their tensor-valued counterparts, typically use a different mapping in order to facilitate enforcing appropriate continuity of only normal or tangential components.

The reference element paradigm is typically employed for elements that satisfy a kind of \emph{equivalence}, where the reference element basis functions map directly to the physical element basis functions under the coordinate change.  However, for many elements (both classical and modern), this is not the case, as the degrees of freedom (DOFs) are `mixed up' by the pullback in a cell-dependent manner.  In~\cite{kirby-zany}, Kirby developed a general theory for obtaining the proper basis in the case of affinely mapped cells.  This approach gives the correct nodal basis as a linear combination of the mapped reference element basis functions.  This linear combination turns out to be sparse, meaning that applying the theory incurs only a small additional cost during the finite element computation, ensuring that the use of these more exotic elements composes neatly with existing high-level automatic code generation software. Kirby and Mitchell~\cite{finat-zany} generalized the Firedrake~\cite{Rathgeber:2016} code stack to generate and employ this transformation, giving results for Morley, Hermite, Argyris, and Bell elements. The basis on the reference cell for each of these elements is constructed using FIAT~\cite{Kirby:2004}.  

In this paper, we extend the transformation theory presented in~\cite{kirby-zany} to Piola-mapped $\hdiv$ elements, both vector- and tensor-valued. Specifically, we focus on the Mardal--Tai--Winther~\cite{mardal2002} and Arnold--Winther~\cite{arnold2002,arnold2003} elements. The vector-valued Mardal--Tai--Winther element provides a discretization of $\hdiv$ which is nonconforming in $H^1$, and works robustly for the scale of operators that interpolate between Stokes and Darcy flow. Thanks to its discrete Korn inequality, it serves as a locking-free elasticity element~\cite{mardal2006observation}. It also provides a low-order divergence-free discretization of the incompressible Stokes and Navier--Stokes equations~\cite{John2017}. It has been implemented for a small number of standalone numerical experiments~\cite{danisch2007}, 
and implemented partially in Diffpack~\citep[Ch.~4]{langtangen2012}
and SyFi~\citep[Ch.~15]{fenics2012_book}
(a symbolic alternative to FIAT) and envisioned as part of FIAT~\citep[Ch.~3]{fenics2012_book}. 
These latter implementations were not fully usable on arbitrary unstructured meshes, due to the lack of transformation theory which we now provide in this work. 

The tensor-valued Arnold--Winther elements are used to discretize the Cauchy stress tensor in the stress-displacement formulation of planar linear elasticity. After decades of effort by many researchers, these elements solved a fundamental challenge 
of numerical elasticity in providing a stable and convergent discretization of the Hellinger--Reissner principle, which moreover enforces symmetry of the stress tensor exactly. Despite being conceived almost two decades ago, the 
prior implementations of these Arnold--Winther elements have either required the explicit element-by-element construction of the basis~\cite{carstensenAWimp}, or again, have been for one-off numerical experiments, such as for stress reconstruction or the design of error indicators~\cite{Goriely2013,gedicke2018,nicaise2008,Riedlbeck2017}. Again, FIAT has had partial support for this element, but the lack of reference mappings has prevented its use in practice~\cite[Ch.~3]{fenics2012_book}.

We begin in Section~\ref{sec:piola} with notation used to define geometry and Piola mappings, and summarize some results on how normal and tangential components are transformed. 
In Section~\ref{sec:vector}, we survey the applications which motivate the use of the Mardal--Tai--Winther element, before describing several other vector-valued elements of interest;
analogously in Section~\ref{sec:tensor}, we briefly review the role of the Cauchy stress tensor in linear elasticity and its approximation with finite elements, before defining the Arnold--Winther stress elements.
We survey the transformation theory of Kirby~\cite{kirby-zany} and show how it applies in the context of Piola elements in Section~\ref{sec:transform}, which is the main contribution of the paper.
A discussion of the application of the Arnold--Winther elements to linear elasticity follows in Section~\ref{sec:discretize_H-R}; in particular, we prove convergence of Nitsche's method for the imposition of traction boundary conditions.
Section~\ref{sec:FEEC} offers an interpretation of Piola transformations in terms of the finite element exterior calculus (FEEC), and reviews the role of FEEC in the development of multigrid preconditioners, as discussed in
Section~\ref{sec:smoothers}, in which we apply patch-based multigrid algorithms to precondition the canonical Stokes--Darcy and Hellinger--Reissner systems for the Mardal--Tai--Winther and Arnold--Winther elements respectively. 
We briefly discuss the implementation within the Firedrake code stack and present some numerical examples for our newly-enabled elements in Section~\ref{sec:numerics}.

\vspace{-3mm}\section{Piola transformation theory}\label{sec:piola}
\vspace{-2mm}

\subsection{Notation and preliminaries}
\vspace{-2mm}

Let $\mathbb{M} = \mathbb{R}^{d\times d}$ denote the space of $d\times d$ matrices, and $\mathbb{S} = \mathbb{R}^{d\times d}_{\text{sym}}$ its symmetric subspace. We shall employ the standard notation for the Sobolev space $H^k(\Omega;X)$ (or $L^2(\Omega;X)$ when $k = 0$), with domain $\Omega\subset\mathbb{R}^d$ and codomain $X$, and associated norm $\|\cdot\|_{k,\Omega}$ or simply $\|\cdot\|_k$, as well as the standard spaces $\hdiv, H(\Div;\mathbb{M}), H(\Div;\mathbb{S}), H(\rot), H(\rot;\mathbb{S}), \hcurl$, and $H(\curl;\mathbb{S})$,
where the divergence, rotation, or curl 
of a tensor field is defined row-wise. 
Differential operators with subscript $h$, such as $\nabla_h$, are meant element-wise.
The symbol $\lesssim$ denotes inequality up to a constant which may depend on mesh regularity but not mesh spacing $h$.

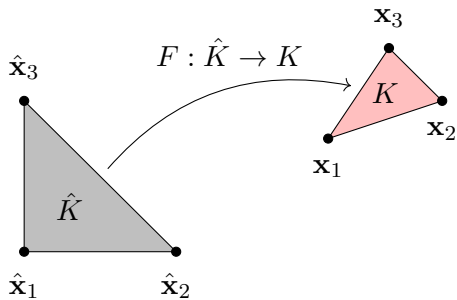
\begin{figure}[H]
  \vspace{-3mm}
  \centering
  \begin{tikzpicture}
    \draw[fill=lightgray] (0,0) coordinate (xhat1)
    -- (2,0) coordinate(xhat2)
    -- (0,2) coordinate(xhat3)--cycle;
    \foreach \pt\labpos\lab in {xhat1/below/\hat{\mathbf{x}}_1, xhat2/below/\hat{\mathbf{x}}_2, xhat3/above/\hat{\mathbf{x}}_3}{
      \filldraw (\pt) circle(.6mm) node[\labpos=1.5mm, fill=white]{$\lab$};
    }
    \draw[fill=pink] (4.0, 1.5) coordinate (x1)
    -- (5.5, 2.0) coordinate (x2)
    -- (4.8, 2.7) coordinate (x3) -- cycle;
    \foreach \pt\labpos\lab in {x1/below/\mathbf{x}_1, x2/below/\mathbf{x}
_2, x3/above/\mathbf{x}_3}{
      \filldraw (\pt) circle(.6mm) node[\labpos=1.5mm, fill=white]{$\lab$};
    }
    \draw[->] (1.1, 1.1) to[bend left] (4.3, 2.2);
    \node at (2.7, 2.65) {$F:\hat{K}\rightarrow K$};
    \node at (0.6,0.6) {$\hat{K}$};
    \node at (4.75, 2.1) {$K$};
  \end{tikzpicture}
    \vspace{-4mm}
    \caption{Affine mapping from a reference cell \(\hat{K}\) to a typical cell \( K \).}\label{fig:affmap}
    \vspace{-5mm}
\end{figure}

Let $\hat{K}\subset\mathbb{R}^2$ be a reference triangle with vertices
$\{\hat{\mathbf{x}}_i\}_{i=1}^{3}$ and let $K$ be a typical element
with vertices $\{\mathbf{x}_i\}_{i=1}^{3}$.  We assume the diffeomorphism
$F: \hat{K} \rightarrow K$ to be affine, as shown in 
Fig.~\ref{fig:affmap}.  We number the edges of a triangle by the convention in~\cite{rognes2009efficient} as follows.
Edge $i$ of any triangle excludes the vertex $i$ and is oriented from the lower-numbered vertex to the higher-numbered one.  Each edge $\be$ is a column vector running between two vertices of a triangle, and is also used to denote the edge as a set.  Its magnitude is the edge length, and we define the unit tangent to edge $\be$ as $\bt = \tfrac{\be}{\| \be \|}$.  To each edge, we associate its normal $\bn = R \bt$, where
\vspace{-1mm}\begin{equation}
    R = \begin{bmatrix}  & 1 \\ -1 &  \end{bmatrix}
\vspace{-1mm}\end{equation}
is clockwise rotation. Quantities and differential operators with $\widehat{~~}$ are defined analogously for the reference element. Note that the normal to a triangle may point inward or outward according to this convention, and we allow $F$ to reverse orientation, but triangles sharing an edge will automatically agree on the direction (so that local and global direction coincide).  This greatly simplifies the required logic for assembly of finite element DOFs involving normal components.  This convention can be extended to simplices of any dimension~\cite{arnold2009geometric,rognes2009efficient}. 

Given any function $\hat{\phi}$ defined on $\hat{K}$, we define its \emph{pullback} by
\vspace{-1mm}\begin{equation}\label{eq:pullback}
    \phi \coloneqq F^{-*}(\hat{\phi}) \coloneqq \hat{\phi} \circ F^{-1}.
\vspace{-1mm}\end{equation}
The pullback can be used to evaluate  $\phi$ at a point $\mathbf{x} = F\left(\hat{\mathbf{x}}\right)$ by
\vspace{-1mm}\begin{equation}
  \phi(\mathbf{x}) = \hat{\phi}(\hat{\mathbf{x}}).
\vspace{-1mm}\end{equation}
For scalar-valued spaces, the pullback provides an isomorphism between $H^m(\hat{K})$ and $H^m(K)$ for $m\in\mathbb{Z}_{\geq 0}$;
its use is very natural in the context of affine equivalent
finite elements such as Lagrange or Crouzeix--Raviart.  One can define its basis functions on a fixed reference element $\hat{K}$ and obtain the basis functions on any $K$ by pullback. This can be adapted, with some technicalities, when affine equivalence fails~\cite{kirby-zany}.

When working with finite elements discretizing Sobolev spaces such as $\hdiv$, it is frequently helpful to work with the contravariant Piola mapping.  Piecewise polynomial vector fields will lie in $\hdiv$ if and only if they have continuous normal components across interfaces~\citep[Prop.~2.1.2]{BBF}; applying this row-wise gives an analogous statement for $H(\Div;\mathbb{M})$ and $H(\Div;\mathbb{S})$. Hence, finite elements for $\hdiv$-based spaces tend to include normal components as degrees of freedom, and these will be preserved by using the appropriate mapping.

\vspace{-2mm}\begin{definition}
    Let $\hat{K}, K \subset \mathbb{R}^d$, let $F:\hat{K}\to K$ be a diffeomorphism, and let $J(\bhx) = \hat{\nabla}F(\hat{\bx})$ be its Jacobian.  The \emph{contravariant Piola map} takes $\hat{\Phi}:\hat{K}\to\mathbb{R}^d$ to $\Phi:K\to\mathbb{R}^d$ defined by
    \vspace{-1mm}\begin{equation}\label{eq:single_Piola}
        \Fdiv(\hat{\Phi}) \coloneqq \tfrac{1}{\det J} J\hat{\Phi} \circ F^{-1}.
    \vspace{-1mm}\end{equation}
    The \emph{(double) contravariant Piola map} takes $\htau:\hat{K}\to\mathbb{S}$ to $\tau:K\to\mathbb{S}$ by
    \vspace{-1mm}\begin{equation}\label{eq:double_con_Piola}
        \Fdivdiv(\hat{\tau}) \coloneqq \tfrac{1}{\left( \det J\right)^2} J\left( \hat{\tau} \circ F^{-1} \right) J^T.
    \vspace{-1mm}\end{equation}
\end{definition}

\noindent The Piola transforms $\Fdiv:H(\Div,\hat{K})\to H(\Div, K), \Fdivdiv:H(\Div,\hat{K};\mathbb{S})\to H(\Div,K;\mathbb{S})$ are 
well-defined and isomorphisms. Analogous considerations hold for $H(\rot), H(\rot;\mathbb{S}), \hcurl, H(\curl;\mathbb{S})$, tangential traces, and the covariant Piola mappings:
\vspace{-2mm}\begin{definition}
    The \emph{covariant Piola map} takes $\hat{\Phi}:\hat{K}\to\mathbb{R}^d$ to $\Phi:K\to\mathbb{R}^d$ by
    \vspace{-1mm}\begin{equation}\label{eq:covariant}
        \mathcal{F}^{\curl}(\hat{\Phi}) \coloneqq J^{-T}\hat{\Phi}\circ F^{-1}.
    \vspace{-1mm}\end{equation}
    The \emph{(double) covariant Piola map} takes $\htau:\hat{K}\to\mathbb{S}$ to $\tau:K\to\mathbb{S}$ by
    \vspace{-1mm}\begin{equation}\label{eq:double_cov_Piola}
        \mathcal{F}^{\curl,\curl}(\htau) \coloneqq J^{-T}\left(\hat{\tau}\circ F^{-1}\right) J^{-1}.
    \vspace{-1mm}\end{equation}
\end{definition}

\noindent Unlike in the usual definition of the Piola map, the absolute value of the determinant is not taken in~\eqref{eq:single_Piola}; this is explained in~\cite{rognes2009efficient}, which gives a thorough exposition of transforming and assembling vector finite elements in the case that the spaces and degrees of freedom are exactly preserved under the Piola mappings. 
The central contribution of this paper is to adapt the techniques of~\cite{kirby-zany} to handle $\hdiv$ elements that are \emph{not} exactly preserved under contravariant Piola mapping.  This occurs, for example, when elements use degrees of freedom involving edge tangents or vertex values in addition to edge/face normals and internal ones.  This enables us to develop general and robust software implementations of these elements.

\vspace{-3mm}\subsection{Transformation of components and their moments}\label{sec:transformed_components}
\vspace{-2mm}

Now, we review the relationship between the reference and physical element normals and tangents when $F$ is affine.  Suppose edge $\be = \bx_i - \bx_j$ connects vertices $i$ and $j$ of the physical element.  Since $F(\hat{\bx}_i) = \bx_i$, we have
\vspace{-1mm}\begin{equation}
  \be = \bx_i - \bx_j = J\left( \hat{\bx}_i - \hat{\bx}_j \right) = J \hat{\be},
\vspace{-1mm}\end{equation}
and equivalently,
\vspace{-1mm}\begin{equation}
  \label{eq:mapt}
  \| \be \| \bt =  \| \hat{\be} \| J \hat{\bt}.
\end{equation}
In a similar way, we use the identity $\tfrac{1}{\det J} R J R^T = J^{-T}$ for all $2 \times 2$ matrices to obtain
\vspace{-1mm}\begin{equation}
\label{eq:mapn}
\| \be \| \bn = R \be = R J \hat{\be} = \left(\det J\right) J^{-T} R \hat{\be} = \left( \det J \right) J^{-T} \| \hat{\be} \| \hat{\bn}.
\vspace{-1mm}\end{equation}
In other words,~\eqref{eq:mapt} shows that the Jacobian maps reference tangents onto physical ones (up to scale factor), but a different mapping given in~\eqref{eq:mapn} is required to obtain the physical normal from the reference one. 

These calculations directly inform how the Piola mapping works.  Let us consider some $\hat{\Phi}$ defined over $\hat{K}$ and its image under the Piola map $\Phi = \Fdiv(\hat{\Phi})$.  At each point $\bx$ on edge $\be \subset \partial K$ the Piola mapping preserves the normal component scaled by edge length, for we have
\vspace{-1mm}\begin{equation}\label{eq:piolan}
  \| \be \| \Phi(\bx) \cdot \bn
   = \left( \tfrac{1}{\det J} J \hat{\Phi}(\hat{\bx}) \right)^T
  \left( \| \hat{\be} \|  \left( \det J \right) J^{-T} \hat{\bn} \right) 
  = \| \hat{\be} \| \hat{\Phi}(\hat{\bx}) \cdot \hat{\bn}.
\vspace{-1mm}\end{equation}
This calculation also specifies how integral moments of normal components against some function $\mu(\bx) = \hat{\mu}(F^{-1}(\bx))$ along edges behave under Piola mapping.  Changing variables from $\be$ to $\hat{\be}$ introduces a Jacobian factor of $\| \be \| / \| \hat{\be} \|$, so that
\vspace{-1mm}\begin{equation}\label{eq:piolanmom}
    \int_\be \left( \Phi(s) \cdot \bn \right) \mu(s)\ds
    = \int_{\hat{\be}} \tfrac{\| \hat{\be} \|}{\|\be\|} \left( \hat{\Phi}(\sh) \cdot \hat{\bn} \right) \hat{\mu}(\sh) \tfrac{\| \be \|} {\| \hat{\be} \|}\dshat
    = \int_{\hat{\be}} \left( \hat{\Phi}(\sh)\cdot \hat{\bn} \right) \hat{\mu}(\sh) \dshat.
\vspace{-1mm}\end{equation}
On the other hand, the Piola map fails to preserve the tangential direction:
\vspace{-1mm}\begin{equation}\label{eq:piolat}
    \| \be \| \Phi(\bx) \cdot \bt
    = \left( \tfrac{1}{\det J} J \hat{\Phi}(\hat{\bx}) \right)^T \left( \| \hat{\be} \| J \hat{\bt} \right)
    = \tfrac{\| \hat{\be} \|}{\det J} \hat{\Phi}(\hat{\bx})^T J^T J \hat{\bt}.
\vspace{-1mm}\end{equation}
It is useful to further resolve $J^T J \hat{\bt}$ and hence $\|\be\| \Phi(\bx) \cdot \bt$ in terms of the reference normal and tangential components.  We define
\vspace{-1mm}\begin{equation}\label{eq:ghat}
    \hat{G} = \begin{bmatrix} \hat{\bn} & \hat{\bt} \end{bmatrix},
\vspace{-1mm}\end{equation}
which is an orthogonal matrix.  Using $\mathbb{I} = \hat{G} \hat{G}^T$ gives
\vspace{-1mm}\begin{equation}
    J^T J \hat{\bt} 
    = \hat{G} \left( \hat{G}^T J^T J \hat{\bt} \right)
    = \tilde{\alpha} \hat{\bn} + \tilde{\beta} \hat{\bt},
\vspace{-1mm}\end{equation}
where
\vspace{-1mm}\begin{equation}
    \begin{bmatrix}
        \tilde{\alpha} \\ \tilde{\beta}
    \end{bmatrix}
    = \hat{G}^T J^T J \hat{\bt}.
\vspace{-1mm}\end{equation}
Now, with $\alpha = \tfrac{\tilde{\alpha}}{\det J}$ and
$\beta = \tfrac{\tilde{\beta}}{\det J}$, we have
\vspace{-1mm}\begin{equation}
  \| \be \| \Phi(\bx) \cdot \bt =
  \| \hat{\be} \| \hat{\Phi}(\hat{\bx}) \cdot \left( \alpha \hat{\bn} + \beta \hat{\bt} \right).
\vspace{-1mm}\end{equation}
So, tangential moments are not preserved under the Piola map, but we have
\vspace{-1mm}\begin{equation}\label{eq:piolatmom}
    \int_\be \left( \Phi(s)\cdot \bt \right) \mu(s) \ds
    = \int_{\hat{\be}} \left( \hat{\Phi}(\sh) \cdot \left( \alpha \hat{\bn} + \beta \hat{\bt} \right) \right) \hat{\mu}(\sh) \dshat.
\vspace{-1mm}\end{equation}

These properties show that vector-valued elements with boundary degrees of freedom involving only normal components can be directly mapped by the Piola transform.  Such elements include the Raviart--Thomas~\cite{Raviart1977}, Brezzi--Douglas--Marini~\cite{brezzi1985two}, and Brezzi--Douglas--Fortin--Marini~\cite{brezzi1987efficient} elements.  On the other hand, the Mardal--Tai--Winther element includes tangential as well as normal boundary degrees of freedom, and we can expect complications in mapping them.

We can also use the calculations~\eqref{eq:mapt} and~\eqref{eq:mapn} to determine what happens to the various components of a tensor under the Piola map.  For some symmetric tensor-valued function $\hat{\tau}: \hat{K} \rightarrow \mathbb{S}$ and $\tau = \Fdivdiv(\hat{\tau})$, the double-normal component is preserved in an analogous sense to the vector case:
\vspace{-1mm}\begin{equation}\label{eq:piolann}
\begin{split}
    \| \be \|^2 \bn^T \tau(\bx) \bn  =  \left( \| \hat{\be}  \| \left( \det J \right) J^{-T} \hat{\bn} \right)^T
    \left( \tfrac{1}{\left( \det J \right)^2} J \hat{\tau}(\hat{\bx}) J^T \right)
    \left( \| \hat{\be} \| \left( \det J \right) J^{-T} \hat{\bn} \right) 
    =  \| \hat{\be} \|^2 \hat{\bn}^T \hat{\tau}(\hat{\bx}) \hat{\bn}.
\end{split}
\vspace{-1mm}\end{equation}
On the other hand, the normal/tangential and double tangential components are not preserved, as direct calculations verify:
\vspace{-1mm}\begin{equation}\label{eq:piolant}
\| \be \|^2 \bn^T \tau(\bx) \bt  \\
= \frac{\| \hat{\be} \|^2}{\det J} \hat{\bn}^T \hat{\tau}(\hat{\bx}) \left( J^T J \right) \hat{\bt},
\vspace{-1mm}\end{equation}
and
\vspace{-1mm}\begin{equation}
    \| \be \|^2 \bt^T \tau(\bx) \bt =
    \frac{\|\hat{\be}\|^2}{\left( \det J \right)^2}
    \hat{\bt}^T \left( J^T J \right) \hat{\tau}(\hat{\bx}) \left( J^T J \right) \hat{\bt}.
\vspace{-1mm}\end{equation}
The tangential-normal component $\| \be \|^2 \bt^T \tau(\bx) \bn$ follows from~\eqref{eq:piolant} by symmetry.
Moreover, we can keep $\hat{G}$ as in~\eqref{eq:ghat} and further resolve these in terms of the reference normal/tangential directions by
\vspace{-1mm}\begin{equation}
\begin{split}
\| \be \|^2 \bn^T \tau(\bx) \bt
 = \|\hat{\be} \|^2 \hat{\bn}^T \hat{\tau}(\hat{\bx}) \left( \alpha \hat{\bn} + \beta \hat{\bt} \right) 
 = \| \hat{\be} \|^2 \left( \alpha \hat{\bn}^T \hat{\tau}(\hat{\bx}) \hat{\bn}
+ \beta \hat{\bn}^T \hat{\tau}(\hat{\bx}) \hat{\bt} \right),
\end{split}
\vspace{-1mm}\end{equation}
and
\vspace{-1mm}\begin{equation}
\| \be \|^2 \bt^T \tau(\bx) \bt
= 
\| \hat{\be} \|^2 \left( \alpha^2 \hat{\bn}^T \hat{\tau}(\hat{\bx}) \hat{\bn}
+ 2 \alpha \beta \hat{\bn}^T \hat{\tau}(\hat{\bx}) \hat{\bt}
+ \beta^2 \hat{\bt}^T \hat{\tau}(\hat{\bx}) \hat{\bt} \right).
\vspace{-1mm}\end{equation}

\vspace{-1mm}\begin{remark}
    The tangential and tangential-tangential 
    traces are well-defined at least on polynomial subspaces, and are included here for completeness.
\vspace{-1mm}\end{remark}

Like in the vector-valued case, we can also show how moments transform.  Suppose $\hat{\tau} \in L^2(\hat{K}; \mathbb{S})$ with $\tau = \Fdivdiv(\hat{\tau})$ and $\hat{\mu} \in L^2(\hat{\be})$. Integrating both sides of~\eqref{eq:piolann} and accounting for the Jacobian $\ds = \tfrac{\|\be\|}{\|\hat{\be}\|}\dshat$ gives
\vspace{-1mm}\begin{equation}
  \label{eq:piolannmom}
\| \be \| \int_{\be} \left( \bn^T \tau(s) \bn \right) \mu \ds
= 
\| \hat{\be} \| \int_{\hat{\be}}\left( \hat\bn^T \hat\tau(\sh) \hat\bn \right) \hat\mu \dshat.
\vspace{-1mm}\end{equation}
Similarly, 
\vspace{-1mm}\begin{equation}
  \label{eq:piolantmom}
\| \be \| \int_{\be} \left( \bn^T \tau(s) \bt \right) \mu \ds
= 
\|\hat{\be} \| \left( \alpha \int_{\hat{\be}}\left( \hat\bn^T \hat\tau(\sh) \hat\bn \right) \hat\mu \dshat
+ \beta \int_{\hat{\be}}\left( \hat\bn^T \hat\tau(\sh) \hat\bt \right) \hat\mu \dshat \right),
\vspace{-1mm}\end{equation}
and
\vspace{-1mm}\begin{equation}
    \label{eq:piolattmom}
\begin{split}
\| \be \| \int_{\be} \left( \bt^T \tau(s) \bt \right) \mu \ds
= \| \hat{\be} \| 
 \left(
    \alpha^2 \int_{\hat{\be}}\left( \hat\bn^T \hat\tau(\sh) \hat\bn \right) \hat\mu \dshat \right. & \left. + 2 \alpha \beta \int_{\hat{\be}}\left( \hat\bn^T \hat\tau(\sh) \hat\bt \right) \hat\mu \dshat  \right. \\
& \left. +\beta^2 \int_{\hat{\be}}\left( \hat\bt^T \hat\tau(\sh) \hat\bt \right) \hat\mu \dshat \right).
\end{split}
\vspace{-1mm}\end{equation}

It follows that tensor-valued elements with boundary DOFs involving only
normal-normal components can be directly mapped by the Piola transform; this
is for example true of the Tangential-Displacement Normal-Normal Stress element~\cite{TDNNS1} and the Hellan--Herrmann--Johnson element~\cite{Hellan1967, Herrmann1967, Johnson1973}. Since the Arnold--Winther elements
also have boundary DOFs of other kinds, application of the Piola map alone will not result in the correct dual basis on each physical cell.

\vspace{-2mm}\begin{remark}
    We briefly comment on the covariant case for $\hcurl$ elements. As is well-known, in 2D, $H(\rot)$ and $\hdiv$ are isomorphic by simple rotation, although the spaces are quite different in 3D; indeed, in dimensions strictly greater than 2, the analogous identification is always impossible for the vector-valued spaces~\citep[Prop.~4.1]{lithesis}. It follows that the transformation theory for $H(\rot)$ and $\hcurl$ elements, at least in 2D, may be similar \textit{mutatis mutandis} to what we present in this paper.
\vspace{-2mm}\end{remark}

\vspace{-3mm}\section{Vector-valued elements}\label{sec:vector}
\vspace{-2mm}

\subsection{Motivations}\label{sec:stokes--darcy}
\vspace{-2mm}

The Mardal--Tai--Winther (MTW) element~\cite{mardal2002} provides a uniformly robust element for both $\hdiv$ and (nonconforming) $H^1$. It is a single discretization that is stable and accurate for both Stokes and Darcy flow as well as the scale of operators interpolating between them, including systems arising from transient Stokes-like flow with small timesteps. Due to its lack of vertex degrees of freedom, it is also appropriate for the interface conditions arising in coupled Stokes--Darcy flow~\cite{karper2009unified}. In the context of elasticity, it admits a discrete Korn inequality~\cite{mardal2006observation} and serves as a locking-free element for primal linear elasticity~\cite{kuchta2019singular}, and for Biot's consolidation model in poroelasticity~\cite{lee2018robust}.

It is especially significant for the purposes of incompressible flow that the MTW space $V_h$ may be paired with the space $Q_h$ of piecewise constants, so that the resulting velocity-pressure pair $V_h\times Q_h$ satisfies $\Div V_h\subset Q_h$, a highly desirable property since it implies mass conservation (i.e.~exactly divergence-free velocities) and pressure robustness~\cite{John2017}.

The MTW element offers a nonconforming discretization of the two-dimensional Stokes complex, a smoother variant of the classical $L^2$ de Rham complex~\cite{tai2006}:
\vspace{-1mm}\begin{equation}\label{eq:stokes_complex}
    \begin{tikzcd}
        0 \arrow{r} & H^2(\Omega)/\mathbb{R} \arrow{r}{\curl} & H^1(\Omega;\mathbb{R}^2) \arrow{r}{\Div} & L^2(\Omega) \arrow{r} & 0.
    \end{tikzcd}
\vspace{-1mm}\end{equation}
\noindent While it does not form a true subcomplex due to its nonconforming nature, this relates the MTW element to elements for other spaces and gives stability and convergence theory for mixed problems; indeed, the original motivation for the construction of the MTW space was that its divergence-free subspace is exactly the $\curl$ of a related $H^2$-nonconforming element, studied in~\cite{nilssen2001robust} and which we denote by $W_h$:
\vspace{-1mm}\begin{equation}\label{eq:MTW_complex}
    \begin{tikzcd}
        0 \arrow{r} & W_h/\mathbb{R} \arrow{r}{\curl} & V_h \arrow{r}{\Div} & Q_h \arrow{r} & 0.
    \end{tikzcd}
\vspace{-1mm}\end{equation}
\noindent Both sequences~\eqref{eq:stokes_complex},~\eqref{eq:MTW_complex} are exact whenever $\Omega$ is simply connected. Additionally, this perspective guides the construction of robust smoothers for optimal-order multilevel algorithms~\cite{FKMW}, as we shall discuss in Sections~\ref{sec:FEEC} and~\ref{sec:smoothers}.

\vspace{-3mm}\subsection{The Mardal--Tai--Winther element}
\vspace{-2mm}

In this Section, we define several other vector-valued $\hdiv$ elements, to discuss their transformations: the Raviart--Thomas~\cite{Raviart1977}, Brezzi--Douglas--Marini~\cite{brezzi1985two}, as well as the Mardal--Tai--Winther~\cite{mardal2002} elements. Schematics of these elements are shown in Figure~\ref{fig:vecels}. 

\begin{figure}[H]
    \vspace{-3mm}
    \begin{subfigure}[t]{0.32\textwidth}
        \centering
        \begin{tikzpicture}[scale=3.0]
            \draw[line width=0.5mm,fill=olive] (0,0) -- (1, 0) -- (0.5, 0.866) -- cycle;
            \foreach \i/\j/\n/\t in {0.5/0.0/0.0/-1,
                                 0.75/0.433/0.866/0.5,
                                 0.25/0.433/-0.866/0.5}{
             \draw[line width=0.4mm,-Latex] (\i, \j) -- (\i+\n/4, \j+\t/4);
                                }
        \end{tikzpicture}
        \vspace{-2mm}
        \caption{Raviart--Thomas~\cite{Raviart1977}}
  \end{subfigure}
  \hfill
        \begin{subfigure}[t]{0.32\textwidth}
        \centering
        \begin{tikzpicture}[scale=3.0]
            \draw[line width=0.5mm,fill=cyan] (0,0) -- (1, 0) -- (0.5, 0.866) -- cycle;
            \foreach \i/\j/\n/\t in {0.6667/0.0/0.0/-1,
                                 0.3333/0.0/0.0/-1,
                                 0.8333/0.2887/0.866/0.5,
                                 0.6667/0.5774/0.866/0.5,
                                 0.1667/0.2887/-0.866/0.5,
                                 0.3333/0.5774/-0.866/0.5}{
             \draw[line width=0.4mm,-Latex] (\i, \j) -- (\i+\n/4, \j+\t/4);
                                }
        \end{tikzpicture}
        \vspace{-2mm}
            \caption{Brezzi--Douglas--Marini~\cite{brezzi1985two}}
  \end{subfigure}
  \hfill
        \begin{subfigure}[t]{0.32\textwidth}
        \centering
        \begin{tikzpicture}[scale=3.0]
            \draw[line width=0.5mm,fill=pink] (0,0) -- (1, 0) -- (0.5, 0.866) -- cycle;
            \foreach \i/\j/\n/\t in {0.6667/0.0/0.0/-1,
                                 0.3333/0.0/0.0/-1,
                                 0.8333/0.2887/0.866/0.5,
                                 0.6667/0.5774/0.866/0.5,
                                 0.1667/0.2887/-0.866/0.5,
                                 0.3333/0.5774/-0.866/0.5}{
             \draw[line width=0.4mm,-Latex] (\i, \j) -- (\i+\n/4, \j+\t/4);
                                }
             \foreach \x/\y/\X/\Y/\none/\ntwo in {0.3/0.0433/0.7/0.0433/0/-1,
                                     0.8/0.2598/0.6/0.6062/0.866/0.5,
                                     0.4/0.6062/0.2/0.2598/-0.866/0.5}{
                                         \draw[line width=0.4mm,-Latex] (\x-\none/14,\y-\ntwo/14) {} -- (\X-\none/14,\Y-\ntwo/14) {};
                                     }
        \end{tikzpicture}
        \vspace{-2mm}
            \caption{Mardal--Tai--Winther~\cite{mardal2002}}
  \end{subfigure}
  \vspace{-4mm}
  \caption{Some triangular vector $\hdiv$ elements.  Arrows represent moments of a particular vector component (normal or tangential) along a facet.}\label{fig:vecels}
  \vspace{-5mm}
\end{figure}

Raviart--Thomas elements (the original $\hdiv$ element) are based on the function space
\vspace{-1mm}\begin{equation}
    \text{RT}_m(K) = \mathcal{P}_m(K;\mathbb{R}^2) + \bx \mathcal{P}_m(K).
\vspace{-1mm}\end{equation}
This is somehow the smallest space of polynomials over $K$ such that the divergence maps exactly onto $\mathcal{P}_m(K)$. The edge degrees of freedom are taken as moments of the normal components up to degree $m$. The Brezzi--Douglas--Marini elements~\cite{brezzi1985two} and Brezzi--Douglas--Fortin--Marini elements~\cite{brezzi1987efficient} have larger spaces than RT, have similar edge degrees of freedom, and so transform in a very similar way.
All these elements form equivalent families under the contravariant Piola map.

We let $\{ \mu_i \}_{i=0}^m$ be some hierarchically ordered basis for
polynomials of degree $m$ on the unit interval (e.g.~monomials or
Legendre polynomials).  (The hierarchical property will not be important
for the Raviart--Thomas element, but we shall use it for the
Mardal--Tai--Winther space.)  We will assume that $\mu_0 = 1$.  Let $\mu^k_i$ be the mapping of $\mu_i$
from the unit interval to edge $\be_k$ of triangle $K$ by
standard pullback.

The Raviart--Thomas degrees of freedom include normal moments on each edge:
\vspace{-1mm}\begin{equation}\label{eq:normmom}
    \ell^{\bn,i,k}(\Phi) = \int_{\be_k} \left( \Phi \cdot \bn \right) \mu_i^k \ds
\vspace{-1mm}\end{equation}
for each edge $\be_k$.  Now let
$
  \{ \hat{\phi}_i \}_{i=1}^{m(m+1)/2}
$
be some hierarchically ordered basis for $\mathcal{P}_{m-1}(\hat{K})$ (e.g.~the
bivariate monomials or Dubiner basis) and $\phi_i = \hat{\phi}_i \circ F^{-1}$
its pullback.  We then let
$\Phi_{k,i} \in L^2(K; \mathbb{R}^2)$ be the vector-valued function with $\phi_i$ in
component $k$ and 0 in the rest.  Defining this, the rest of the degrees of freedom are
\vspace{-1mm}\begin{equation}\label{eq:internmom}
    \ell^{k,i}(\Phi) = \int_K \Phi \cdot \Phi_{k,i} \dx.
\vspace{-1mm}\end{equation}
Then, taking the union of all of these gives the nodes for the Raviart--Thomas element.
If we construct the nodal basis $\{\hat{\Psi}_i\}$ for $\text{RT}_{m}(\hat{K})$ or $\text{BDM}_m(\hat{K})$ on a reference element, then a suitable basis for $\text{RT}_m$ or $\text{BDM}_m$ on any $K$ is obtained by applying $\Fdiv$ to each member of the space.  This follows directly from the preservation of normal sense~\eqref{eq:mapn}.

In this context, the Mardal--Tai--Winther element~\cite{mardal2002} is an example of an $\hdiv$ element that includes both normal and tangential degrees of freedom and so requires a more complex reference mapping.  On each cell $K$, the MTW space is
\vspace{-1mm}\begin{equation}
    \text{MTW}(K) = \left\{ \Phi \in \mathcal{P}_3(K; \mathbb{R}^2)~\vert~\Div \Phi \in \mathcal{P}_0(K), \left(\Phi \cdot \bn \right)|_{\be_k} \in \mathcal{P}_1(\be_k), 1\leq k\leq 3 \right\},
\vspace{-1mm}\end{equation}
which includes all linear polynomials (and some higher-order terms) and has dimension exactly nine.  Six degrees of freedom are the constant and linear moments of the normal components of each edge, as in~\eqref{eq:normmom}.  Instead of internal moments like~\eqref{eq:internmom} for Raviart--Thomas, the remaining three are the constant tangential moment on each edge $\be_k$:
\vspace{-1mm}\begin{equation}
    \ell^{\bt,0,k}(\Phi) = \int_{\be_k} \Phi \cdot \bt \ds.
\vspace{-1mm}\end{equation}
The globally defined MTW space is therefore fully continuous across edges in the normal component, and continuous in a weaker sense in the tangential component, hence is `almost' $H^1$-conforming; moreover, application to problems posed in $H^1$ are simplified in that no facet integral penalty terms are required.

This element is somewhat more expensive than $\hdiv$ elements with comparable accuracy.  Like $\text{BDM}_1$, it pairs with piecewise constant pressures but has 50\% more degrees of freedom (three per edge rather than two). It gives a second order approximation with only one more degree of freedom than $\text{RT}_2$, and unlike $\text{RT}_2$ serves also as a $H^1$-nonconforming element, but $\text{RT}_2$ pairs with piecewise linear pressures.

\vspace{-3mm}\section{Tensor-valued elements}\label{sec:tensor}
\vspace{-2mm}

\subsection{Elastic stress and the numerical enforcement of symmetry}\label{sec:elastic_stress}
\vspace{-2mm}

For practical applications of linear elasticity, one is often interested in computing the Cauchy stress tensor with at least as much accuracy as the displacement, and hence, in mixed stress-displacement formulations in which the stress is computed directly rather than, for example, via numerical differentiation after the fact.
The mixed dual stress-displacement formulation achieves this, 
and offers a way to alleviate the well-known numerical phenomenon of volumetric locking in the incompressible limit.

Let $\Omega\subset\mathbb{R}^d$ be a polygonal elastic body; we focus on the planar case $d = 2$. Given Lam\'e parameters $\lambda, \mu > 0$, the compliance tensor $\mathcal{A}:\mathbb{S}\to\mathbb{S}$ is, in the homogeneous isotropic case, defined as
\vspace{-1mm}\begin{equation}
    \mathcal{A}\tau \coloneqq \frac{1}{2\mu}\left(\tau - \frac{\lambda}{2\mu + d\lambda}(\tr\tau)\mathbb{I}\right).
\vspace{-1mm}\end{equation}
Let $\Gamma_D, \Gamma_N$ partition $\partial\Omega$, and let $\langle\cdot,\cdot\rangle_{\Gamma_D}$ denote the $(H^{-1/2}\times H^{1/2}_{00})(\Gamma_D;\mathbb{R}^d)$ dual pairing. Given an external body force $f\in L^2(\Omega;\mathbb{R}^d)$, displacement data $u_0 \in H^{1/2}_{00}(\Gamma_D;\mathbb{R}^d)$, and traction data $g \in H^{-1/2}(\Gamma_N;\mathbb{R}^d)$, we consider the Hellinger--Reissner system, which seeks a stress-displacement pair $(\sigma,u)\in H(\Div;\mathbb{S})\times L^2(\Omega;\mathbb{R}^d)$ satisfying
\vspace{-1mm}\begin{equation}\label{eq:traction_cond}
    \sigma \bn = g \text{ on }\Gamma_N\text{ in the trace sense, }
\vspace{-1mm}\end{equation}
and which are critical points of the Hellinger--Reissner functional
\vspace{-1mm}\begin{equation}\label{eq:H-R_fnl}
    \mathcal{H}(\sigma, u) \coloneqq \int_\Omega \frac{1}{2}\mathcal{A}\sigma:\sigma + (\Div\sigma - f)\cdot u~\dx - \langle\sigma\bn, u_0\rangle_{\Gamma_D}.
\vspace{-1mm}\end{equation}
A stationary point $(\sigma, u)$ satisfies
\vspace{-1mm}\begin{equation}\label{eq:H-R_stat}
    \begin{split}
        \int_\Omega \mathcal{A}\sigma:\tau + (\Div\tau)\cdot u + (\Div\sigma)\cdot v~\dx = \int_\Omega f\cdot v~\dx + \langle\tau\bn, u_0\rangle_{\Gamma_D}\\
        \forall~(\tau, v) \in H(\Div;\mathbb{S})\times L^2(\Omega;\mathbb{R}^d)\text{ with }\tau\bn\vert_{\Gamma_N} = 0,
    \end{split}
\vspace{-1mm}\end{equation}
i.e.~the saddle point system
\vspace{-1mm}\begin{subequations}\label{eq:H-R_strongform}
    \begin{alignat}{2}
        \mathcal{A}\sigma &= \varepsilon(u) \quad &&\text{ in }\Omega, \label{eq:constit} \\
        \Div\sigma &= f \quad &&\text{ in }\Omega, \label{eq:equil}\\
        u &= u_0 \quad &&\text{ on }\Gamma_D, \label{eq:Dircond}\\
        \sigma \bn &= g \quad &&\text{ on }\Gamma_N, \label{eq:Neumcond}
    \end{alignat}
\vspace{-1mm}\end{subequations}
where $\varepsilon$ denotes the symmetric gradient (or linearized strain). 
Here, the constitutive relation~\eqref{eq:constit} is associated with the minimization among stress fields of an appropriate complementary energy, and is formally equivalent to the standard ``Hooke's law" $\sigma = \mathbb{C}\varepsilon(u) = 2\mu\varepsilon(u) + \lambda(\tr\varepsilon(u))\mathbb{I}$ (where $\mathbb{C} = \mathcal{A}^{-1}$ denotes the elasticity tensor), which in primal formulations is typically used to \textit{define} (and thus eliminate) the stress. However, only the former stress-strain relation~\eqref{eq:constit} remains valid in the incompressible limit $\lambda\gg\mu$.

That the Cauchy stress field $\sigma$ is symmetric is equivalent to the conservation of angular momentum, and is highly desirable but notoriously difficult to preserve in finite element discretizations; a further challenge arises from the requirement of $H(\Div;\mathbb{M})$-conformity, which ensures that the traction forces on a mesh face shared between two elements are in equilibrium. 
Before defining the $H(\Div;\mathbb{S})$-discretizing, exactly symmetric Arnold--Winther (AW) elements, we
briefly review other approaches to the numerical enforcement of symmetry; 
for further references, see~\citep[Ch.~9.3--9.4]{BBF}.

Conventionally, symmetry is enforced weakly via $L^2(\Omega;\mathbb{M})$-orthogonality to a skew-symmetric subspace; this approach encapsulates 
the PEERS element~\cite{PEERS}, and the weakly symmetric successors to the AW elements~\cite{AFW2007}. 
Weak imposition of symmetry also fits naturally into the FOSLS (first order system least-squares) method
by penalization of asymmetry~\cite{Cai, Muller}.

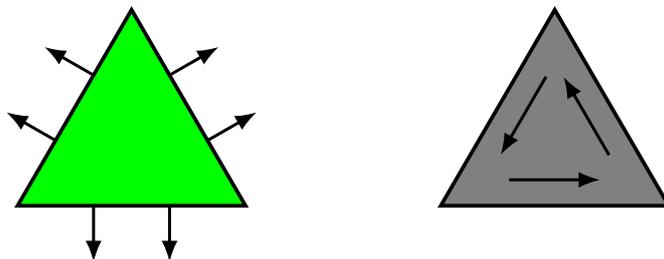
\begin{figure}[H]
    \vspace{-3mm}
    \centering
    \begin{subfigure}{0.32\textwidth}
        \centering
        \begin{tikzpicture}[scale=3.0]
            \draw[line width=0.5mm,fill=green] (0,0) -- (1, 0) -- (0.5, 0.866) -- cycle;
            \foreach \i/\j/\n/\t in {0.6667/0.0/0.0/-1,
                                 0.3333/0.0/0.0/-1,
                                 0.8333/0.2887/0.866/0.5,
                                 0.6667/0.5774/0.866/0.5,
                                 0.1667/0.2887/-0.866/0.5,
                                 0.3333/0.5774/-0.866/0.5}{
             \draw[line width=0.4mm,-Latex] (\i, \j) -- (\i+\n/4, \j+\t/4);
                                }
        \end{tikzpicture}
    \end{subfigure}
    \begin{subfigure}{0.32\textwidth}
        \centering
        \begin{tikzpicture}[scale=3.0]
            \draw[line width=0.5mm,fill=gray] (0,0) -- (1, 0) -- (0.5, 0.866) -- cycle;
            \foreach \x/\y/\X/\Y/\none/\ntwo in {0.3/0.0433/0.7/0.0433/0/-1,
                                     0.8/0.2598/0.6/0.6062/0.866/0.5,
                                     0.4/0.6062/0.2/0.2598/-0.866/0.5}{
                                         \draw[line width=0.4mm,-Latex] (\x-\none/14,\y-\ntwo/14) {} -- (\X-\none/14,\Y-\ntwo/14) {};
                                     }
             \draw[line width=0.4mm,->,color=white] (0.5,-0.02) -- (0.5,-0.25);
        \end{tikzpicture}
    \end{subfigure}
    \vspace{-4mm}
    \caption{The Tangential-Displacement and Normal-Normal-Stress mixed element~\cite{TDNNS1}. Thin arrows on the tensor diagram refer to a single tensor component in the given direction (e.g.~$\bn^T \tau \bn$).}\label{fig:TDNNS_diagram}
    \vspace{-5mm}
\end{figure}

A very significant contribution are the TDNNS elements due to Sch\"oberl and Pechstein (n\'ee Sinwel)~\cite{TDNNS1};
the stress space consists of the symmetric tensor-valued polynomials of degree $k\geq 1$ whose normal-normal component is
continuous across each edge, paired with a displacement space of the same degree whose tangential component is continuous; schematics are provided in Figure~\ref{fig:TDNNS_diagram}.

\begin{figure}[H]
    \vspace{-3mm}
    \centering
    \begin{subfigure}{0.32\textwidth}
        \centering
        \begin{tikzpicture}[scale=3.0]
            \draw[line width=0.5mm,fill=yellow] (0,0) -- (1, 0) -- (0.5, 0.866) -- cycle;
            \foreach \i/\j/\n/\t in {0.5/0.0/0.0/-1,
                                 0.75/0.433/0.866/0.5,
                                 0.25/0.433/-0.866/0.5}{
             \draw[line width=0.4mm,-Latex] (\i, \j) -- (\i+\n/4, \j+\t/4);
                                }
        \end{tikzpicture}
    \end{subfigure}
    \vspace{-4mm}
    \caption{The Hellan--Herrmann--Johnson element~\cite{Hellan1967, Herrmann1967, Johnson1973}.}\label{fig:HHJ_diagram}
    \vspace{-5mm}
\end{figure}
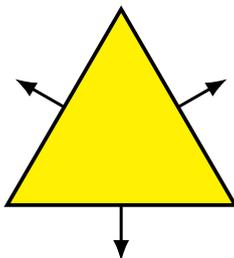

Another exactly symmetric element is the Hellan--Herrmann--Johnson (HHJ) element, 
consisting 
of symmetric matrix-valued
polynomials of degree $k \geq 0$ with continuity of the normal-normal
components (that is, $\bn^T \tau \bn$ is continuous across edges).
A schematic is provided in Figure~\ref{fig:HHJ_diagram}.

\noindent The TDNNS elements are implemented in Netgen/NGSolve, while the HHJ element was introduced into FEniCS in~\cite{hale2018simple} and later into Firedrake;
both implementations were possible without the theory we develop in this paper since both have
boundary degrees of freedom involving only the normal-normal
components and hence, in light of~\eqref{eq:piolann}, can be Piola-mapped in
a standard way. 

\vspace{-3mm}\subsection{The Arnold--Winther elements}\label{sec:AW_elts}
\vspace{-2mm}

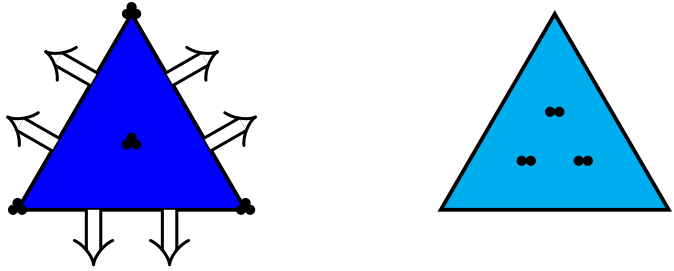
\begin{figure}[H]
    \vspace{-3mm}
    \centering
    \begin{subfigure}[t]{0.32\textwidth}
        \centering
        \begin{tikzpicture}[scale=3.0] 
            \draw[line width=0.5mm,fill=blue] (0,0) -- (1, 0) -- (0.5, 0.866) -- cycle;
            \foreach \i/\j in {0/0,1/0,0.5/0.866}{
                \foreach \dx/\dy in {0.02/0,-0.02/0,0/0.03}{
                    \draw[fill=black] (\i+\dx,\j+\dy) circle (0.02);
                }
            }
            \foreach \i/\j/\n/\t in {0.6667/0.0/0.0/-1,
                                 0.3333/0.0/0.0/-1,
                                 0.8333/0.2887/0.866/0.5,
                                 0.6667/0.5774/0.866/0.5,
                                 0.1667/0.2887/-0.866/0.5,
                                 0.3333/0.5774/-0.866/0.5}{
             \draw[line width=0.4mm, double, double distance=1.5mm,-implies] (\i, \j) -- (\i+\n/4, \j+\t/4);
                                }
            \foreach \dx/\dy in {0.02/0.0,-0.02/0.0,0.0/0.03}{
            \draw[fill=black] (0.5+\dx,0.2887+\dy) circle (0.02);
                                }
        \end{tikzpicture}
    \end{subfigure}
    \begin{subfigure}[t]{0.32\textwidth}
        \centering
        \begin{tikzpicture}[scale=3.0]
            \draw[line width=0.5mm,fill=cyan] (0,0) -- (1, 0) -- (0.5, 0.866) -- cycle;
            \foreach \dx in {0.02,-0.02}{
                \foreach \a/\b in {0.375/0.2165,0.625/0.2165,0.5/0.433}{
                    \draw[fill=black] (\a+\dx,\b) circle (0.02);
                                }
                            }
             \draw[line width=0.4mm, double, double distance=1.5mm, color=white,-implies] (0.5,-0.02) -- (0.5,-0.25);
        \end{tikzpicture}
    \end{subfigure}
    \vspace{-4mm}
    \caption{The conforming Arnold--Winther element~\cite{arnold2002}, and the discontinuous Lagrange element for displacement with which it is paired. Thick arrows refer to both components of a tensor in a given direction (e.g.~$\tau \bn$ or $\tau \bt$). A clutch of circles either indicates internal moments of each unique component or evaluation of each component at a 
    point.}\label{fig:AWc_diagram}
    \vspace{-5mm}
\end{figure}

Arnold and Winther proposed two exactly symmetric elements discretizing $H(\Div;\mathbb{S})$, one conforming~\cite{arnold2002}
and the other nonconforming~\cite{arnold2003}. We first consider the conforming element.
In the lowest-order case, the space is the symmetric matrix-valued quadratic polynomials, augmented by solenoidal cubics:
\vspace{-1mm}\begin{equation}
\begin{aligned}
    \text{AW\textsuperscript{c}}(K) &= \{\tau\in\mathcal{P}_3(K;\mathbb{S})~\vert~\Div\tau\in\mathcal{P}_1(K;\mathbb{R}^2)\} \\
            &= \mathcal{P}_2(K;\mathbb{S}) + \{\tau\in\mathcal{P}_3(K;\mathbb{S})~\vert~\Div\tau = 0\}.
\end{aligned}
\vspace{-1mm}\end{equation}
This space has dimension 24, and the degrees of freedom are:
\vspace{-3mm}\begin{itemize}\setlength\itemsep{-1mm}
    \item the values of each component of $\tau$ at each vertex of $K$,
    \item the moments of degree 0 and 1 of the normal-normal and
      normal-tangential components of $\tau$ on each edge,
    \item the constant moment of each component of $\tau$ over $K$.
\vspace{-3mm}\end{itemize}
The associated local displacement space is $V_h(K) = \mathcal{P}_1(K;\mathbb{R}^2)$ which is 6-dimensional, with DOFs given by, for example, the values of the 2 components at 3 non-collinear points interior to $K$.
Element diagrams are provided in Figure~\ref{fig:AWc_diagram}.

To introduce notation for the stress DOFs, we let
\vspace{-1mm}\begin{equation}
    \ell^{k,i,j}(\tau) = \tau_{i,j}(\bx_k)
\vspace{-1mm}\end{equation}
denote evaluation of a particular component of the input $\tau$ at one
of the vertices of triangle $K$.
Similar to the nodes for vector-valued elements, we define the boundary moments by
\vspace{-1mm}\begin{equation}
    \ell^{\bn,\bs,i, k}(\tau)
    = \int_{\be_k} \left( \bn \cdot \tau \bs \right) \mu_i \ds,
\vspace{-2mm}\end{equation}
where $\bs\in\{\bn, \bt\}$;  for the HHJ element, we thus 
have nodes
$ \{ \ell^{\bn,\bn, i, k} \}_{i=0}^1$ for each edge $\be_k$ of the triangle.

For the interior nodes, 
denoting by $\{\phi_i\}_i$ a hierarchically ordered basis over $K$ with $\phi_0 = 1$, we let
$\tau_{n,i,j}$ denote the tensor with $\phi_n$ in the $i,j$
entry and 0 elsewhere -- a natural injection into the space of tensor-valued orthogonal polynomials.  We let
\vspace{-1mm}\begin{equation}\label{eq:internal_DOF}
    \ell^{n,i,j}(\tau) = \int_{K} \tau : \tau_{n,i,j} \dx.
\vspace{-1mm}\end{equation}
\noindent Thus, DOFs for the AW\textsuperscript{c} space are given in full by $\{\ell^{k, 1,1},\ell^{k, 1, 2},\ell^{k, 2, 2}\}_{k=1}^3, \{\ell^{\bn,\bn,i,k},\ell^{\bn,\bt,i,k}\}_{i=0,k=1}^{1,\hspace{3mm}3}$, $\{\ell^{0,1,1},\ell^{0,1,2},\ell^{0,2,2}\}$.

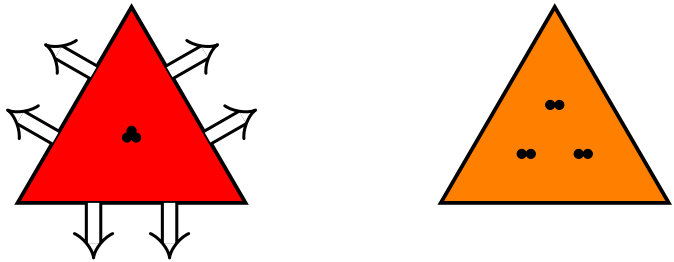
\begin{figure}[H]
    \vspace{-3mm}
    \centering
    \begin{subfigure}{0.32\textwidth}
        \centering
        \begin{tikzpicture}[scale=3.0] 
            \draw[line width=0.5mm,fill=red] (0,0) -- (1, 0) -- (0.5, 0.866) -- cycle;
            \foreach \i/\j/\n/\t in {0.6667/0.0/0.0/-1,
                                 0.3333/0.0/0.0/-1,
                                 0.8333/0.2887/0.866/0.5,
                                 0.6667/0.5774/0.866/0.5,
                                 0.1667/0.2887/-0.866/0.5,
                                 0.3333/0.5774/-0.866/0.5}{
             \draw[line width=0.4mm, double, double distance=1.5mm,-implies] (\i, \j) -- (\i+\n/4, \j+\t/4);
                                }
            \foreach \dx/\dy in {0.02/0.0,-0.02/0.0,0.0/0.03}{
            \draw[fill=black] (0.5+\dx,0.2887+\dy) circle (0.02);
                                }
        \end{tikzpicture}
    \end{subfigure}
    \begin{subfigure}{0.32\textwidth}
        \centering
        \begin{tikzpicture}[scale=3.0]
            \draw[line width=0.5mm,fill=orange] (0,0) -- (1, 0) -- (0.5, 0.866) -- cycle;
            \foreach \dx in {0.02,-0.02}{
                \foreach \a/\b in {0.375/0.2165,0.625/0.2165,0.5/0.433}{
                    \draw[fill=black] (\a+\dx,\b) circle (0.02);
                                }
                            }
             \draw[line width=0.4mm, double, double distance=1.5mm, color=white,-implies] (0.5,-0.02) -- (0.5,-0.25);
        \end{tikzpicture}
    \end{subfigure}
    \vspace{-4mm}
    \caption{The nonconforming Arnold--Winther mixed element~\cite{arnold2003}.}\label{fig:AWnc_diagram}
    \vspace{-5mm}
\end{figure}

The nonconforming stress element introduced in~\cite{arnold2003}
avoids the somewhat unusual feature of vertex degrees of freedom and
gives a cheaper (though slightly less accurate) method based on
the space 
\vspace{-1mm}\begin{equation}
    \text{AW\textsuperscript{nc}}(K) =
    \left\{ \tau \in \mathcal{P}_2(K; \mathbb{S})~\vert~\bn_k \cdot \tau \bn_k \in \mathcal{P}_1(\be_k), \ \ 1 \leq k
    \leq 3 \right\}.
\vspace{-1mm}\end{equation}
This space, consisting of symmetric tensors of quadratic polynomials subject to a degree reduction in the normal-normal component on each edge, is 15-dimensional, and is determined by $\{\ell^{\bn,\bn,i,k},\ell^{\bn,\bt,i,k}\}_{i=0,k=1}^{1,\hspace{3mm}3}$, $\{\ell^{0,1,1},\ell^{0,1,2},\ell^{0,2,2}\}$, i.e.
\vspace{-3mm}\begin{itemize}\setlength\itemsep{-1mm}
    \item the constant and
    linear moments of both the normal-normal
    and normal-tangential components on each edge,
    \item the constant
    moments of the three unique components on the interior.
\vspace{-2mm}\end{itemize}
It is again paired with the displacement space $\mathcal{P}_1(K;\mathbb{R}^2)$; 
element diagrams are provided in Figure~\ref{fig:AWnc_diagram}.

The Arnold--Winther elements were one of the first products of the then embryonic finite element exterior calculus~\cite{AFW2006}. Their novelty is that they form an exact subcomplex of the stress elasticity complex in 2D, with commuting cochain projections:
\vspace{-1mm}\begin{equation}\label{eq:AW_complex}
    \begin{tikzcd}
        0 \arrow{r} & \mathcal{P}_1(\Omega) \arrow{r}{\subset} \arrow{d}{\text{id}} & H^2(\Omega) \arrow{r}{\airy} \arrow{d}{I_h} & H(\Div;\mathbb{S}) \arrow{r}{\Div} \arrow{d}{\Pi_h} & L^2(\Omega;\mathbb{R}^2) \arrow{r} \arrow{d}{P_h} & 0\\
        0 \arrow{r} & \mathcal{P}_1(\Omega) \arrow{r}{\subset} & Q_h \arrow{r}{\airy_h} & \Sigma_h \arrow{r}{\Div_h} & V_h \arrow{r} & 0.
    \end{tikzcd}
\vspace{-1mm}\end{equation}
\noindent Here, if $(\Sigma_h, V_h)$ are either the conforming or nonconforming AW elements, $Q_h$ is either the Argyris space (another element requiring nonstandard transformations) or a certain nonconforming approximation of $H^2(\Omega)$ respectively, $\airy$ denotes the Airy stress operator which we will define in~\eqref{eq:airy_stress}, $I_h, \Pi_h$ are appropriate densely defined interpolants, and $P_h$ is the $L^2$-projection, then it is proved in~\cite{arnold2002, arnold2003} (see also~\cite{Wang_PhD}) that the above diagram is commuting, with exact rows whenever $\Omega$ is simply connected.
It is these structure-preserving properties which are used to prove unisolvency of dual bases, prove error estimates, and clarify links between their elements and $H^2$-(non)conforming elements.
Indeed, under the unifying framework, the elasticity system~\eqref{eq:H-R_strongform} is the mixed Hodge Laplacian boundary value problem for a specific choice of Hilbert complex -- more precisely, a segment of the stress complex such that the base space $L^2(\Omega;\mathbb{S})$ is endowed with the energy norm induced by the compliance tensor $\mathcal{A}$~\citep[p.~107]{Arnold_FEEC}.

\vspace{-2mm}\begin{remark}\label{rem:H-Z}
    Many methods proposed after the seminal papers of Arnold and Winther~\cite{arnold2002, arnold2003} were, in essence, attempts to construct improved versions of the Arnold--Winther elements. Among these, perhaps the most efficient are the Hu--Zhang elements~\cite{Hu, HZ1, HZ2}, which strongly and stably enforce symmetry in any dimension; in the 2D conforming case, they have 18 local DOFs in comparison to AW\textsuperscript{c}'s 24. The same authors recently provided a unified analysis of methods which strongly impose symmetry~\cite{Hong}, using a generalized 4-field formulation of which both the Hu--Zhang and AW elements, among many others, are special cases. These elements also admit an interpretation in terms of the FEEC~\cite{christiansen2018}. The Hu--Zhang elements employ tangential DOFs, and the Piola transformation theory which they therefore require would be very similar to that for the AW elements which we present here.
\vspace{-2mm}\end{remark}

\vspace{-3mm}\section{Transforming Piola-inequivalent elements}\label{sec:transform}
\vspace{-2mm}

\subsection{Reviewing the transformation theory}
\vspace{-2mm}

Equations~\eqref{eq:piolat},~\eqref{eq:piolatmom},~\eqref{eq:piolantmom}, and~\eqref{eq:piolattmom} indicate that finite elements involving tangential degrees of freedom will not map nicely under the Piola transform -- the reference element nodal basis will not map to the physical element one.  In~\cite{kirby-zany}, Kirby developed a theory to address the simpler but analogous situation with standard pullback. The broad strokes of the theory go through essentially unchanged if one replaces the regular pullback $F^{-*}(\hat{\phi}) = \hat{\phi} \circ F^{-1}$ with a Piola pullback; for self-containment, we summarize the key ideas here.

We let $\hat{\mathbf{\Phi}}$ be a (column) vector of elements of a function space $\hat{V}$ over $\hat{K}$, taking values in some space such as $\mathbb{R}^d$ or $\mathbb{S}$, with $\dim V = n$.  Let $F: \hat{K} \rightarrow K$ be the mapping shown in Figure~\ref{fig:affmap} and $\mathcal{F}^*:\hat{V}\to V$ some pullback operator (e.g.~one of the Piola maps) taking $\hat{V}$ into a function space $V$ defined over $K$.  We define $\mathcal{F}^*(\hat{\mathbf{\Phi}})$ to be the componentwise application of the pullback:
\vspace{-1mm}\begin{equation}\label{eq:Fvec}
    \left(\mathcal{F}^*(\hat{\mathbf{\Phi}})\right)_j \coloneqq \mathcal{F}^*(\hat{\mathbf{\Phi}}_j).
\vspace{-1mm}\end{equation}
\noindent Then, if we take $\hat{\mathbf{\Psi}}$
to contain the nodal basis functions
of $\hat{V}$ on the reference element and $\mathbf{\Psi}$
to contain the nodal basis 
of $V$ defined on $K$, an immediate consequence of equations such as~\eqref{eq:piolatmom} or~\eqref{eq:piolantmom} is that 
\vspace{-1mm}\begin{equation}\label{eq:nodal_basis_neq}
    \mathcal{F}^*(\hat{\mathbf{\Psi}}) \neq \mathbf{\Psi},
\vspace{-1mm}\end{equation}
\noindent and note that equality still fails to hold up to permutation or scaling.
If the pullback operator preserves the function space, 
providing an isomorphism between the instantiations $V$ and $\hat{V}$, then
$\mathcal{F}^*(\hat{\mathbf{\Psi}})$ will be a basis -- just not the nodal one.  Consequently, there must exist a (cell-dependent) matrix $M$ such that
\vspace{-1mm}\begin{equation}
    \mathbf{\Psi} = M \mathcal{F}^*(\hat{\mathbf{\Psi}}).
\vspace{-1mm}\end{equation}
In this case, one can simply compute the pullback of the reference basis and multiply by $M$.  If $M$ is quite sparse, this is a small additional cost in the finite element computation.  With some effort, the theory can be extended to handle the situation when the pullback does not send $\hat{V}$ onto $V$, but none of our examples in the present work require this.

However, as we see from subsection~\ref{sec:transformed_components}, it is natural to consider the effect of the pullback on the DOFs, rather than directly on the nodal basis. To consider instead how the dual spaces are transformed, denote by $N = \{\ell_i\}_i$ and $\hat{N} = \{\hat{\ell}_i\}_i$ (row) vectors of functionals on $V$ and $\hat{V}$ respectively, and define the \textit{dual} or \textit{push-forward} $\Fs:V^*\to\hat{V}^*$ associated with $\mathcal{F}^*$ by
\vspace{-1mm}\begin{equation}\label{eq:push-forward}
    \Fs(\ell) \coloneqq \ell\circ\mathcal{F}^*.
\vspace{-1mm}\end{equation}
Now identifying an $n$-vector $\mathbf{\Phi}\in V^n$ with its componentwise image $\mathbf{\Phi}^{**}\in (V^{**})^n$ under the canonical embedding into the bidual, it induces an evaluation operator on vectors of dual functionals via the `outer product' $N\mathbf{\Phi} \coloneqq \mathbf{\Phi}^{**}\otimes N$, where
\vspace{-1mm}\begin{equation}
    (\mathbf{\Phi}^{**}\otimes N)_{i,j} \coloneqq \mathbf{\Phi}_j^{**}(\ell_i) = \ell_i(\mathbf{\Phi}_j).
\vspace{-1mm}\end{equation}
This allows us to extend the push-forward~\eqref{eq:push-forward} to vectors of functionals componentwise either in a manner analogous to~\eqref{eq:Fvec}, or equivalently as
\vspace{-1mm}\begin{equation}
    \Fs(N) \coloneqq N\mathcal{F}^* = \mathcal{F}^*(\cdot)^{**}\otimes N.
\vspace{-1mm}\end{equation}
Choosing $N, \hat{N}$ to be the DOFs, then the Kronecker property by which they are characterized is then simply expressed as
\vspace{-1mm}\begin{equation}
    \Psi^{**}\otimes N = \mathbb{I} = \hat{\Psi}^{**}\otimes\hat{N}.
\vspace{-1mm}\end{equation}
Moreover, this clarifies the transformation of the primal basis: expanding the desired physical basis in terms of the mapped reference basis, we have
\vspace{-1mm}\begin{equation}
    \mathbf{\Psi}_j = \sum_{k=1}^n \theta_k^{(j)}\mathcal{F}^*(\hat{\mathbf{\Psi}}_k), \qquad j = 1,\ldots, n,
\vspace{-1mm}\end{equation}
so it follows from the Kronecker property that
\vspace{-1mm}\begin{equation}\label{eq:vandermonde}
    Q\Theta = \mathbb{I},
\vspace{-1mm}\end{equation}
where
\vspace{-1mm}\begin{equation}
    Q = \mathcal{F}^*(\hat{\mathbf{\Psi}})^{**}\otimes N = \Big(\ell_i(\mathcal{F}^*(\hat{\mathbf{\Psi}}_k))\Big)_{i, k = 1}^n
\vspace{-1mm}\end{equation}
can be interpreted as a generalized Vandermonde matrix and
\vspace{-1mm}\begin{equation}
    \Theta = Q^{-1} = \Big(\theta_k^{(j)}\Big)_{k,j=1}^n
\vspace{-1mm}\end{equation}
has the vector of coefficients for $\mathbf{\Psi}_j$ with respect to $\{\mathcal{F}^*(\hat{\mathbf{\Psi}}_k)\}_{k}$ in column $j$.

Now by~\eqref{eq:nodal_basis_neq} and unisolvency, we have
\vspace{-1mm}\begin{equation}
    \mathcal{F}^*(\hat{\Psi})^{**}\otimes N \neq \Psi^{**}\otimes N = \mathbb{I},
\vspace{-1mm}\end{equation}
so
\vspace{-1mm}\begin{equation}
    \Fs(N)(\hat{\Psi}) = \mathcal{F}^*(\hat{\Psi})^{**}\otimes N \neq \mathbb{I} = \hat{\Psi}^{**}\otimes\hat{N} = \hat{N}\hat{\Psi},
\vspace{-1mm}\end{equation}
so in particular,
\vspace{-1mm}\begin{equation}
    \Fs(N) \neq \hat{N},
\vspace{-1mm}\end{equation}
but as before, there must be an invertible $P$ with
\vspace{-1mm}\begin{equation}
    \hat{N} = P\Fs(N).
\vspace{-1mm}\end{equation}
An important result of~\citep[Theorem 3.1]{kirby-zany} is that
\vspace{-1mm}\begin{equation}
    M = P^T
\vspace{-1mm}\end{equation}
(as matrices of numbers). That is, it is sufficient to relate the physical nodes and their push-forwards.  As we will see in our examples, it can be more natural to find $P^{-1}$ and then invert it by hand. 

\vspace{-3mm}\subsection{Piola-mapped vector elements}
\vspace{-2mm}

The general approach of the transformation theory developed in~\cite{kirby-zany} is to build the nodal basis on $K$ out of a linear combination of the mapped reference element basis functions. 
As explained in the previous subsection, for us the transformation of basis functions works by duality, considering instead the push-forward of nodes, which follow readily from the calculations~\eqref{eq:piolanmom} and~\eqref{eq:piolatmom}.

We will use block notation in defining our vectors of degrees of freedom.
Let
\vspace{-1mm}\begin{equation}
    N^k = \begin{bmatrix} \ell^{\bn,0,k} & \ell^{\bt,0,k} & \ell^{\bn,1,k} \end{bmatrix}^T
\vspace{-1mm}\end{equation}
contain the three MTW degrees of freedom associated with edge $\be_k$, and then $N$ is block vector with three parts containing the nodes for each edge:
\vspace{-1mm}\begin{equation}
  N = 
  \begin{bmatrix}
      N^1 &
      N^2 &
      N^3
  \end{bmatrix}^T.
\vspace{-1mm}\end{equation}
Similarly, we let $\hat{N}^k$ and $\hat{N}$ contain the corresponding nodes on the reference element. Now, we need to compute the Piola push-forward $\Fs$ of these nodes in terms of the reference element nodes, which effectively builds $P^{-1}$ for the transformation theory. 

Let $\hat{\Phi}$ be defined over $\hat{K}$, $1 \leq k \leq 3$, and $i=0$ or 1.  Then~\eqref{eq:piolanmom} gives
\vspace{-1mm}\begin{equation}
    \begin{split}
        \mathcal{F}_*(\ell^{\bn,i,k})(\hat{\Phi})
         = 
        \ell^{\bn,i,k}(\Fdiv(\hat{\Phi})) 
         = \int_{\be_k} \left( \Fdiv(\hat{\Phi}) \cdot \bn \right) \mu_i \ds 
         = \int_{\hat{\be}_k} \left( \hat{\Phi} \cdot \hat{\bn} \right) \hat{\mu}_i \dshat = \hat{\ell}^{\bn,i,k}(\hat{\Phi}).
    \end{split}
\vspace{-1mm}\end{equation}
So, the normal moment nodes are in fact pushed forward to the corresponding reference element normal moment nodes.  For tangential moments,~\eqref{eq:piolatmom} can be rewritten as
\vspace{-1mm}\begin{equation}
\begin{split}
\mathcal{F}_*(\ell^{\bt,i,k})(\hat{\Phi})
& = \alpha^k
\hat{\ell}^{\bn,i,k}(\hat{\Phi})
+ \beta^k
\hat{\ell}^{\bt,i,k}(\hat{\Phi}),
\end{split}
\vspace{-1mm}\end{equation}
where we have labeled the $\alpha$, $\beta$ with a superscript $k$
indicating they correspond to edge $k$.  We define
\vspace{-1mm}\begin{equation}
W^k =
\begin{bmatrix}
  1 &  &  \\
  \alpha^k & \beta^k &  \\
   &  & 1
\end{bmatrix},
\vspace{-1mm}\end{equation}
and this gives
\vspace{-1mm}\begin{equation}
  \mathcal{F}_*(N^k)
  =
  W^k
  \hat{N}^k.
\vspace{-1mm}\end{equation}
Using this calculation for each edge, we have
\vspace{-1mm}\begin{equation}
  \mathcal{F}_*(N)
  = \begin{bmatrix} W^1 &  &  \\  & W^2 &  \\  &  &
    W^3 \end{bmatrix}
  \hat{N},
\vspace{-1mm}\end{equation}
and the block matrix in this last equation is $P^{-1}$.  It is easy to
invert each $W^k$:
\vspace{-1mm}\begin{equation}
  (W^k)^{-1} = \begin{bmatrix} 1 &  & 
      \\ -\tfrac{\alpha^k}{\beta^k} & \tfrac{1}{\beta^k} &  \\
     &  & 1 \end{bmatrix}
\vspace{-1mm}\end{equation}
and hence obtain $P$ by
\vspace{-1mm}\begin{equation}
  P = \begin{bmatrix} (W^1)^{-1} &  &  \\  & (W^2)^{-1} &  \\  &
     & (W^3)^{-1} \end{bmatrix}.
\vspace{-1mm}\end{equation}

\vspace{-3mm}\subsection{Piola-mapped tensor elements}
\vspace{-2mm}

Applying the 
abstract transformation theory to the Arnold--Winther
stress elements introduced above follows a similar pattern, using
the results in Section~\ref{sec:piola}.
We begin with the nonconforming element, and subsequently consider the inclusion of the vertex DOFs for the conforming element.
The existing implementations of these elements~\cite{carstensenAWimp, carstensen2011} require a separate construction of the basis functions for each element, by inverting the Vandermonde matrix arising from the Piola pullback of the primal basis, as in equation~\eqref{eq:vandermonde}.

For the nonconforming element, we define the vector of nodes associated with edge $\be_k$ by
\vspace{-1mm}\begin{equation}
  N^k =
\begin{bmatrix}
    \ell^{\bn,\bn,0,k} & \ell^{\bn,\bt,0,k} & \ell^{\bn,\bn,1,k} & \ell^{\bn,\bt,1,k} \end{bmatrix}^T.
\vspace{-1mm}\end{equation}
That is, we store the normal and tangential moments of order 0,
followed by those of order 1. We also collect the three
internal degrees of freedom in the vector
\vspace{-1mm}\begin{equation}\label{eq:awinternal}
    N^4 =
    \begin{bmatrix}
        \ell^{0,1,1} &
        \ell^{0,1,2} &
        \ell^{0,2,2}
    \end{bmatrix}^T,
\vspace{-1mm}\end{equation}
so that the nonconforming element has degrees of freedom
\vspace{-1mm}\begin{equation}\label{eq:awncdofs}
    N = \begin{bmatrix} N^1 & N^2 & N^3 & N^4 \end{bmatrix}^T,
\vspace{-1mm}\end{equation}
with a similar definition and ordering of reference element degrees
of freedom $\hat{N}^k$ and $\hat{N}$.

Much as with the MTW element, we can use~\eqref{eq:piolannmom}
and~\eqref{eq:piolantmom} to write
\vspace{-1mm}\begin{equation}
  \mathcal{F}_*( N^k ) = W^k \hat{N}^k,
\vspace{-1mm}\end{equation}
where we now have
\vspace{-1mm}\begin{equation}
  \label{eq:wk}
    W^k = 
    \frac{\|\bhe_k\|}{\|\be_k\|}
    \begin{bmatrix}
        1 &  &  &  \\
        \alpha^k & \beta^k &  &  \\
         &  & 1 &  \\
         &  & \alpha^k & \beta^k
    \end{bmatrix},
\vspace{-1mm}\end{equation}
and we readily invert $W^k$ by
\vspace{-1mm}\begin{equation}
    \left(W^k\right)^{-1} = 
    \frac{\|\be_k\|}{\|\bhe_k\|}
    \begin{bmatrix}
    1 &  &  &  \\
    -\tfrac{\alpha^k}{\beta^k} & \tfrac{1}{\beta_k} &  &  \\
     &  & 1 &  \\
     &  & -\tfrac{\alpha^k}{\beta^k} & \tfrac{1}{\beta_k} \\
  \end{bmatrix}
\vspace{-1mm}\end{equation}
and have that
\vspace{-1mm}\begin{equation}
  P = \begin{bmatrix} (W^1)^{-1} &  &  \\
     & (W^2)^{-1} &  \\
     &  & (W^3)^{-1} \end{bmatrix}
\end{equation}
in a manner very analogous to the MTW element.

Let us apply the push-forward to an internal degree of freedom, recalling that $\phi_0 = 1$:
\vspace{-1mm}\begin{equation}\label{eq:pfint}
  \begin{split}
      \mathcal{F}_*(\ell^{0, i, j})(\hat{\tau}) =
      \int_K \Fdivdiv(\hat{\tau}) : \tau_{0, i, j} \dx 
       = \int_{\hat{K}} \tfrac{1}{\det J} \left( J \hat{\tau} J^T \right)_{ij} \dxhat.
  \end{split}
\vspace{-1mm}\end{equation}
If we expand out the quantity
$
    J \hat{\tau} J^T,
$
we obtain
\vspace{-1mm}\begin{equation}
  \label{eq:jtaujt}
\begin{bmatrix}
  J_{11}^2 \hat{\tau}_{11} + 2 J_{11} J_{12} \hat{\tau}_{12} + J_{12}^2 \hat{\tau}_{22}
  & J_{11} J_{21} \hat{\tau}_{11} + \left( J_{11} J_{22} + J_{12}J_{21}
  \right) \hat{\tau}_{12} + J_{12}J_{22} \hat{\tau}_{22} \\
  * & J_{21}^2 \hat{\tau}_{11} + 2 J_{21}J_{22} \hat{\tau}_{12} + J_{22}^2
  \hat{\tau}_{22}
\end{bmatrix},
\vspace{-1mm}\end{equation}
where the asterisk indicates equality due to symmetry.  Using this in~\eqref{eq:pfint} gives that
\vspace{-1mm}\begin{equation}
  \begin{split}
    \mathcal{F}_*(\ell^{0, 1, 1}) & =
    \tfrac{1}{\det J} \left(
    J_{11}^2 \hat{\ell}^{0, 1, 1}
    + 2 J_{11} J_{12} \hat{\ell}^{0,1,2}
    + J_{12}^2 \hat{\ell}^{0,2,2} \right) \\
    \mathcal{F}_*(\ell^{0, 1, 2}) & = \tfrac{1}{\det J} \left(
    J_{11} J_{21} \hat{\ell}^{0,1,1}
    + \left(J_{11} J_{22} + J_{12} J_{21} \right) \hat{\ell}^{0,1,2}
    + J_{12} J_{22} \hat{\ell}^{0,2,2}
    \right)\\
    \mathcal{F}_*(\ell^{0, 2, 2}) & = \tfrac{1}{\det J} \left(
    J_{21}^2 \hat{\ell}^{0,1,1}
    + 2 J_{21} J_{22} \hat{\ell}^{0,1,2}
    + J_{22}^{2} \hat{\ell}^{0,2,2} \right),
  \end{split}
\vspace{-1mm}\end{equation}
so that with
\vspace{-1mm}\begin{equation}
  \label{eq:wtilde}
  \tilde{W} = 
  \begin{bmatrix}
    J_{11}^2 & 2 J_{11} J_{12} & J_{12}^2 \\
    J_{11}J_{21} & J_{11}J_{22} + J_{12}J_{21} & J_{12} J_{22} \\
    J_{21}^2 & 2 J_{21} J_{22} & J_{22}^2 
  \end{bmatrix}
\vspace{-1mm}\end{equation}
and
$\check{W} = \tfrac{1}{\det J}\tilde{W}$,
we write
\vspace{-1mm}\begin{equation}
  \mathcal{F}_*(N^4) =
  \check{W} \hat{N}^4.
\vspace{-1mm}\end{equation}
This gives the result
\vspace{-1mm}\begin{equation}\label{eq:AWnc_Pinv}
    P^{-1} =
    \begin{bmatrix}
        W^1 &  &  &  \\
         & W^2 &  &  \\
         &  & W^3 &  \\
         &  &  & \check{W}
    \end{bmatrix}
\vspace{-1mm}\end{equation}
which, again, is readily inverted blockwise to find $P$.
As a remark, the $\check{W}$ block is dense, but 
denoting its dependence on $J$ by $\check{W} = \check{W}_J$, it is then easily checked
that $\check{W}_J^{-1} = \check{W}_{J^{-1}}$, i.e.~the inverse can be found by reversing the roles of $K$ and $\hat{K}$.

We now turn to the full conforming Arnold--Winther element.  It has the
same edge and internal nodes to the nonconforming element just considered,
but also includes vertex values. However, these transform in a very similar fashion to the internal moments just discussed.
We collect the pointwise degrees of freedom for each vertex $k$ into a small vector:
\vspace{-1mm}\begin{equation}
    N^{0,k} =
    \begin{bmatrix} \ell^{k,1,1} & \ell^{k,1,2} & \ell^{k,2,2} \end{bmatrix}^T.
\vspace{-1mm}\end{equation}
We also collect the edge degrees of freedom together in the same way as for the
nonconforming element:
\vspace{-1mm}\begin{equation}
    N^{1,k} =
    \begin{bmatrix}
        \ell^{\bn,\bn,0,k} & \ell^{\bn,\bt,0,k} & \ell^{\bn,\bn,1,k} & \ell^{\bn,\bt,1,k}
    \end{bmatrix}^T.
\vspace{-1mm}\end{equation}
Note that we have added extra superscripts of 0 and 1 to indicate the
topological dimension associated with the degrees of freedom.  We also
include the internal degrees of freedom in a separate vector
\vspace{-1mm}\begin{equation}
    N^{2} = \begin{bmatrix} \ell^{0,1,1} & \ell^{0,1,2} & \ell^{0,2,2} \end{bmatrix}^T,
\vspace{-1mm}\end{equation}
and hence order the degrees of freedom by
\vspace{-1mm}\begin{equation}\label{eq:awcdofs}
    N = \begin{bmatrix}
            N^{0,1} & N^{0,2} & N^{0,3} & N^{1,1} & N^{1,2} & N^{1,3} & N^{2}
        \end{bmatrix}^T.
\vspace{-1mm}\end{equation}

The edge and internal degrees of freedom are handled in exactly the same way as for the previous elements. The Piola push-forward of the vertex functional gives 
\vspace{-1mm}\begin{equation}
  \label{eq:vertex}
  \mathcal{F}_*(\ell^{k,i,j})(\hat{\tau}) =
  \ell^{k,i,j}(\Fdivdiv(\hat{\tau}))
    = \left( \tfrac{1}{(\det J)^2} J \hat{\tau}(\hat{\bx}_k) J^T \right)_{ij}.
\vspace{-1mm}\end{equation}
This suggests that (except for very special geometry), the Piola
transform will not send vertex-oriented basis functions on the
reference element to their physical counterpart.  However, using the expansion of $J \hat{\tau} J^T$ given in~\eqref{eq:jtaujt}, a very similar calculation as for the internal moments shows that~\eqref{eq:vertex} gives for each vertex
\vspace{-1mm}\begin{equation}
    \mathcal{F}_*(N^{0,k}) = \breve{W} \hat{N}^{0,k},
\vspace{-1mm}\end{equation}
where now $\breve{W} = \tfrac{1}{(\det J)^2} \tilde{W}$,
for which again $\breve{W}_J^{-1} = \breve{W}_{J^{-1}}$.
Combining all of the push-forwards gives that
\vspace{-1mm}\begin{equation}\label{eq:AWc_Pinv}
  \mathcal{F}_*(N) =
  \begin{bmatrix}
    \breve{W} &  &  &  &  &  &  \\
     & \breve{W} &  &  &  &  &  \\
     &  & \breve{W} &  &  &  &  \\
     &  &  & W^{1,1} &  &  &  \\
     &  &  &  & W^{1,2} &  &  \\
     &  &  &  &  & W^{1,3} &  \\
     &  &  &  &  &  & \check{W}
  \end{bmatrix}
  \hat{N}.
\vspace{-1mm}\end{equation}
Again, the block-diagonal structure of the matrix makes it
straightforward to invert and find $P$.

\vspace{-3mm}\subsection{Scale-invariance and conditioning}\label{sec:scale-inv}
\vspace{-2mm}

If global basis functions are of disparate size as the mesh is refined, this introduces an additional source of ill-conditioning into finite element systems.  This phenomenon occurs for Hermite and certain other affinely-mapped elements~\cite{kirby-zany}, and similar conditions apply to our Piola-mapped elements as well.  The solution is to globally re-scale degrees of freedom so that basis functions are of comparable size.

We  illustrate these issues for the AW\textsuperscript{c} element.
Let $h$ denote the typical diameter of a cell. 
Consider a basis function dual to a vertex DOF; it is unity at one vertex in one component and zero at other vertices, hence $\mathcal{O}(1)$ over the cell.  A basis function dual to an internal DOF integrates to unity over the cell, hence must be of size $\mathcal{O}(h^{-2})$.
Between these cases, a normal moment basis function over an edge must have size $\mathcal{O}(h^{-1})$ over that edge.  Integrating the inner products which populate the global mass matrix, we obtain some entries (corresponding to pairs of vertex basis functions) of size $\mathcal{O}(h^{-2})$, others (corresponding to pairs of internal basis functions) of size $\mathcal{O}(h^2)$, and others in between.  This leads to an $\mathcal{O}(h^{-4})$ condition number of the mass matrix rather than an $\mathcal{O}(1)$ condition number resulting from, say, standard Lagrange polynomials.

This issue can be readily addressed in a \emph{post hoc} fashion by rescaling the global DOFs according to the local mesh size. For example, one can multiply the edge degrees of freedom by $h$ and the internal degrees of freedom by $h^2$, so long as cells sharing a DOF agree on the precise value used in the scaling.
This choice was made in our implementation, but alternatively, one could include some power of facet measure in the definitions of degrees of freedom for the finite element.  Such scale-invariance is implicit in the treatment of AW\textsuperscript{c} by Carstensen et al.~\cite{carstensenAWimp}, the Morley--Wang--Xu elements, and related $H^3$-nonconforming elements due to Wu and Xu~\cite{Xu2019} and would give a slightly different transformation matrix.

\vspace{-3mm}\subsection{Preservation of the constraints}\label{sec:constraints}
\vspace{-2mm}

Our three spaces of primary interest can be characterized as the intersection of the kernels of several linear functionals. 
Note that the constraint on the normal component used to define the MTW space is preserved by the Piola map, due to~\eqref{eq:piolan}. Similarly, the constraint on the normal-normal component of the nonconforming AW\textsuperscript{nc} space is preserved in light of~\eqref{eq:piolann}. The MTW element is moreover constrained by a loss of degree in its divergence;
denoting by $\Phi = \Fdiv(\hat{\Phi})$ the vector Piola map applied to some $\hat{\Phi}$ defined on the reference element, by a direct computation we have
$\Div\Phi = \frac{1}{\det J}\widehat{\Div}~\hat{\Phi}$.
Similarly, the conforming AW\textsuperscript{c} element is also subject to a loss-of-degree constraint on its (now vector-valued) divergence. Denoting $\tau = \Fdivdiv(\hat{\tau})$ for some symmetric $\hat{\tau}$ defined on the reference element, then
$\Div\tau = \frac{1}{(\det J)^2}J\widehat{\Div}~\hat{\tau}$.
It follows that the Piola maps preserve the degree of the divergence in both cases.

These constraints on the divergence, as well as constraints on the degree of the normal components, may be expressed by requiring integral moments against certain orthogonal polynomials to vanish. This is described on the reference element in some detail for Arnold--Winther elements in~\cite{Kirby:2004}.  All these functionals are preserved under the Piola push-forward.  Hence, the physical element function space is actually constructed by pullback, even if the basis functions are not preserved.  This need not be the case -- for example, the constraints on the $C^1$ Bell element are not preserved under affine pullback, requiring a fuller version of the transformation theory~\cite{kirby-zany}.

\vspace{-3mm}\section{Traction conditions in the Hellinger--Reissner principle with Nitsche's method}\label{sec:discretize_H-R}
\vspace{-2mm}

The Arnold--Winther discretization of the Hellinger--Reissner problem~\eqref{eq:H-R_strongform}, in the pure displacement case $\vert\Gamma_N\vert = 0$, for displacement data $u_0$, has discrete weak form given by seeking $(\sigma_h,u_h)$ in one of the AW pairs $\Sigma_h\times V_h$ such that
\vspace{-1mm}\begin{equation}\label{eq:MMS_H-R}
    \begin{aligned}
        \int_\Omega \mathcal{A}\sigma_h:\tau_h + (\Div_h\tau_h)\cdot u_h~\dx &= \int_{\partial\Omega} (\tau_h\bn)\cdot u_0~\ds &\forall~\tau_h\in\Sigma_h, \\
        \int_\Omega (\Div_h\sigma_h)\cdot v_h~\dx &= \int_\Omega f\cdot v_h~\dx &\forall~v_h\in V_h.
    \end{aligned}
\vspace{-1mm}\end{equation}

We now turn to the enforcement of traction conditions~\eqref{eq:traction_cond} for the Hellinger--Reissner problem, which for the dual formulation~\eqref{eq:H-R_stat} are essential (i.e.~strongly enforced); we consider a mixed boundary condition ($0 < |\Gamma_D|,\vert\Gamma_N\vert < |\partial\Omega|$), although the method may be extended to the pure traction case $|\Gamma_D| = 0$ subject to a quotient of the displacement space by the rigid motions and a suitable compatibility condition.

The original AW papers~\cite{arnold2002, arnold2003} only treated the case of the elastic body clamped everywhere on the boundary, i.e.~$\vert\Gamma_N\vert = 0$ with homogeneous displacement data $u_0\equiv 0$. Carstensen et al.~\cite{carstensenAWimp} identified (in particular, inhomogeneous) traction conditions as a substantial practical difficulty associated with the AW elements, due to delicate interdependence between DOFs at the boundary, which moreover depend on the shape of the boundary at a given boundary vertex; they were able to enforce them using nodal interpolation~\cite{carstensen2016}, Lagrange multipliers~\cite{carstensenAWimp}, or elimination of the boundary DOFs by condensation~\cite{carstensen2011}. 
When the boundary condition is mixed but the traction data $g$ is zero, the relevant trial and test space for the stress is 
\vspace{-1mm}\begin{equation}
    H_{0,\Gamma_N}(\Div,\Omega;\mathbb{S}) \coloneqq \{\tau\in H(\Div,\Omega;\mathbb{S})~\vert~\langle\tau\bn,v\rangle_{\Gamma_N} = 0~\forall~v\in H^1_{0,\Gamma_D}(\Omega;\mathbb{R}^2)\}
\vspace{-1mm}\end{equation}
\citep[Remark 2.1.3]{BBF}. Wang~\citep[Lemma III.1]{Wang_PhD} constructed an interpolant of Scott--Zhang type which preserves the traction-free condition on $\Gamma_N$ under an elliptic regularity assumption, and which was used to prove a discrete inf-sup condition for this case, but 
no details of the discrete enforcement of this condition were offered in Pasciak and Wang's application of AW\textsuperscript{c} to the homogeneous pure traction problem in~\cite{pasciak2006}.

There is no clear way to enforce the traction condition in either the Arnold--Winther spaces or the weak formulation of Hellinger--Reissner. 
We therefore advocate a simpler approach, employing the classical Nitsche's method~\cite{Nitsche} to weakly enforce the condition. 
There is little literature on the application of Nitsche's method to the enforcement of essential boundary conditions in dual mixed problems; our approach is 
similar to~\cite{Burman2021, Stenberg2011} for the mixed Poisson problem.

Given traction data $g$, we augment the Hellinger--Reissner functional~\eqref{eq:H-R_fnl} over the discrete spaces $\Sigma_h\times V_h$ with a term incorporating the traction condition (to ensure consistency, since we do not impose the condition on the spaces), and a consistent, quadratic penalty term, seeking the critical point $(\sigma_h,u_h)\in\Sigma_h\times V_h$ of
\vspace{-1mm}\begin{equation}\label{eq:nitsche_HR}
    \mathcal{H}_{h,\gamma}(\sigma_h, u_h) \coloneqq \mathcal{H}(\sigma_h,u_h) -~\int_{\Gamma_N}\br_h\cdot u_h~\ds + \frac{\gamma}{2}\delta_h(\br_h),
\vspace{-1mm}\end{equation}
where 
$h$ denotes the characteristic mesh size, 
$\gamma > 0$ is an $h$-independent penalty parameter, 
$\br_h = \sigma_h\bn - g$ denotes the traction residual, and
\vspace{-1mm}\begin{equation}\label{eq:delta}
    \delta_h(\br_h) \coloneqq \frac{1}{h}\int_{\Gamma_N} \|\br_h\|^2~\ds
\vspace{-1mm}\end{equation}
may be interpreted as a least-squares term penalizing deviation from the traction condition.

The exact traction satisfies $\sigma\bn = g$ in $H^{-1/2}(\Gamma_N;\mathbb{R}^2)$, but a penalty term $\delta_h$ in terms of the dual norm in $H^{-1/2}(\Gamma_N;\mathbb{R}^2)$ would not aid the analysis, nor is it practical to equivalently penalize the Riesz representative of $\br_h$ in $H^{1/2}_{00}(\Gamma_N;\mathbb{R}^2)$, or to work with the linearization of such penalizations. 
To ensure consistency of the $L^2(\Gamma_N;\mathbb{R}^2)$-penalization~\eqref{eq:delta}, 
we assume full elliptic regularity of the stress field by assuming that the solution $(\sigma, u)$ to~\eqref{eq:H-R_stat} satisfies 
$(\sigma, u)\in H^1(\Omega;\mathbb{S})\times H^2(\Omega;\mathbb{R}^2)$,
and 
that $g\in L^2(\Gamma_N;\mathbb{R}^2)$.
The latter assumption typically holds in practice for the traction data (or some discrete approximation thereof), while the former
holds in the homogeneous isotropic case if $\Omega$ is a convex polygon and $\mathcal{A}$ is smooth.

\vspace{-1mm}\begin{remark}
    The exact displacement $u$ formally satisfies an unmixed primal formulation of linear elasticity for which $u\in H^1(\Omega;\mathbb{R}^2)$, but 
    there is no gain in regularity for the stress field by passing to an alternative formulation.
\vspace{-2mm}\end{remark}

\vspace{-2mm}\begin{remark}
    The natural choice of exponent for $h$ in~\eqref{eq:delta} is $+1$ and not $-1$, by dimensionalization and since we expect $\|\br_h\|^2_{0,\Gamma_N}$ to converge one order slower than $\|\br_h\|^2_{-1/2,\Gamma_N}$. 
    However, the term with negative exponent is more naturally interpretable as a penalization, 
    and was found to be more effective in preliminary numerical experiments.
    This also informs the choice of discrete norms~\eqref{eq:mesh-dep_norms} below.
\end{remark}
\vspace{-1mm}
Convergence will be proved only for the AW\textsuperscript{c} element, but a computational example will be included also for AW\textsuperscript{nc}.

Let $\mathcal{T}_h$ denote a 
quasi-uniform
triangulation of $\Omega$, $E_h$ the set of its internal edges, and $\Lambda_h$ the set of its vertices.\footnote{We assume the edges of $\mathcal{T}_h$ align with the partition $\Gamma_D\cup\Gamma_N$, so that each external edge lies entirely in one of these.}
Linearizing the augmented functional~\eqref{eq:nitsche_HR} over the AW\textsuperscript{c} pair, we seek $(\sigma_h, u_h)\in\Sigma_h\times V_h$ satisfying
\vspace{-1mm}\begin{equation}\label{eq:H-R_nitsche}
    \begin{aligned}
        &a_{h,\gamma}(\sigma_h, \tau_h) + b(\tau_h, u_h) \hspace{-3mm}&&= \int_{\Gamma_D}\tau_h\bn\cdot u_0~\ds + \frac{\gamma}{h}\int_{\Gamma_N}g\cdot\tau_h\bn~\ds~&\forall~\tau_h\in\Sigma_h,\\
        &b(\sigma_h, v_h)                              &&= \int_\Omega f\cdot v_h~\dx - \int_{\Gamma_N}g\cdot v_h\ds~&\forall~v_h\in V_h,
    \end{aligned}
\vspace{-1mm}\end{equation}
where
\vspace{-1mm}\begin{align}
    a_{h,\gamma}(\sigma_h,\tau_h) \coloneqq \int_\Omega\mathcal{A}\sigma_h:\tau_h~\dx + \frac{\gamma}{h}\int_{\Gamma_N}\sigma_h\bn\cdot\tau_h\bn~\ds, \quad
    b(\tau_h, v_h) \coloneqq \int_\Omega (\Div\tau_h)\cdot v_h~\dx - \int_{\Gamma_N}\tau_h\bn\cdot v_h~\ds.
\vspace{-1mm}\end{align}
Define the 
$H^1(\Omega;\mathbb{S})\times H^1(\Omega;\mathbb{R}^2)$-based mesh-dependent norms
\vspace{-1mm}\begin{equation}\label{eq:mesh-dep_norms}
    \verth{\tau_h}^2 \coloneqq \|\tau_h\|^2_{0,\Omega} + \frac{1}{h}\|\tau_h\bn\|_{0,\Gamma_N}^2,
    \qquad \|v_h\|_h^2 \coloneqq \|\varepsilon_h(v_h)\|_0^2 + \frac{1}{h}\sum_{\be\in E_h}\|\jump{v_h}\|_{0,\be}^2 + \frac{1}{h}\|v_h\|_{0,\Gamma_D}^2,
\vspace{-1mm}\end{equation}
where $\jump{\cdot}$ on an interior edge $\be = \partial K\cap\partial K'$ denotes the jump $\jump{v_h} \coloneqq v_h\vert_K - v_h\vert_{K'}$, and on an exterior edge $\be\subset\partial\Omega$ denotes the identity.
We now prove 
well-posedness of the augmented discrete formulation~\eqref{eq:H-R_nitsche} using the standard Brezzi conditions~\citep[Section 4.2.3]{BBF}. We take $\gamma\geq 1$.

\vspace{-2mm}\begin{prop}
    The augmented Nitsche system~\eqref{eq:H-R_nitsche} is well-posed, uniformly in $h$, with respect to the norms~\eqref{eq:mesh-dep_norms}.
\vspace{-3mm}\end{prop}
\begin{proof}
There exists $C_{\mathcal{A}} > 0$ with $\int_\Omega\mathcal{A}\tau:\tau~\dx\leq C_{\mathcal{A}}^2\|\tau\|_{0,\Omega}^2~\forall~\tau\in L^2(\Omega;\mathbb{S})$.
For $\sigma_h,\tau_h\in\Sigma_h, v_h\in V_h$, we have
\vspace{-1mm}\begin{equation}
    \vert a_{h,\gamma}(\sigma_h,\tau_h)\vert \leq\max\{C_{\mathcal{A}}^2, \gamma\}\verth{\sigma_h}\verth{\tau_h},
\end{equation}
and
\begin{equation}\label{eq:b_expanded}
    \begin{split}
        b(\tau_h,v_h) &= -\int_{\Gamma_N}\tau_h\bn\cdot v_h~\ds + \sum_{K\in\mathcal{T}_h}\int_\Omega(\Div\tau_h)\cdot v_h~\dx \\
                      &= \sum_K\left(\int_{\partial K}\tau_h\bn\cdot v_h~\ds - \int_K\varepsilon(v_h):\tau_h~\dx\right) - \int_{\Gamma_N}\tau_h\bn\cdot v_h~\ds\\
                      &= -\int_\Omega\varepsilon_h(v_h):\tau_h~\dx + \sum_{\be\in E_h}\int_\be\tau_h\bn\cdot\jump{v_h}\ds + \int_{\Gamma_D}\tau_h\bn\cdot v_h~\ds,
    \end{split}
\end{equation}
so
\begin{equation}
    |b(\tau_h,v_h)| \leq \|\varepsilon_h(v_h)\|_0\|\tau_h\|_0 + \sum_{\be\in E_h}\|\tau_h\bn\|_{0,\be}\|\jump{v_h}\|_{0,\be} + \|\tau_h\bn\|_{0,\Gamma_D}\|v_h\|_{0,\Gamma_D}.
\vspace{-1mm}\end{equation}
By the 
scaling 
$\|\tau_h\bn\|_{0,\be}\lesssim h^{-\frac{1}{2}}\|\tau_h\|_{0,K}~\forall~\be\subset\partial K$, we obtain
\begin{equation}
    |b(\tau_h,v_h)| \lesssim \verth{\tau_h}\|v_h\|_h.
\end{equation}
Clearly, $a_{h,\gamma}$ is $\verth{\cdot}$-coercive on all of $\Sigma_h$ uniformly in $h$, provided $\gamma\geq 1$.
Now fix $0\neq u_h\in V_h$.
Define $\tau_h\in\Sigma_h$ by the AW\textsuperscript{c} DOFs
\begin{subequations}
        \begin{align}
            \tau_h(\bx)&= 0 ~\hspace{-15mm}&&\forall~\bx\in\Lambda_h, \label{eq:tau_vertex}\\
            \int_\be\tau_h\bn\cdot w_h~\ds &= \frac{1}{h}\int_\be\jump{u_h}\cdot w_h~\ds~\hspace{-15mm}&&\forall~\be\in E_h,\forall~\be\subset\Gamma_D, w_h\in\mathcal{P}_1(\be;\mathbb{R}^2),\label{eq:tau_Eh}\\
            \int_\be\tau_h\bn\cdot w_h~\ds &= 0~\hspace{-15mm}&&\forall~\be\subset\Gamma_N, w_h\in\mathcal{P}_1(\be;\mathbb{R}^2),\label{eq:tau_GammaN}\\
            \int_K\tau_h:q_h~\dx &= - \int_K\varepsilon(u_h):q_h~\dx~\hspace{-15mm}&&\forall~K\in\mathcal{T}_h, q_h\in\mathcal{P}_0(K;\mathbb{S}) = \mathbb{S}.\label{eq:tau_internal}
        \vspace{-1mm}\end{align}
\end{subequations}
Choosing $w_h = \jump{u_h}\in\mathcal{P}_1(\be;\mathbb{R}^2)$ in~\eqref{eq:tau_Eh}, 
$q_h = \varepsilon(u_h)\in\mathbb{S}$ in~\eqref{eq:tau_internal} gives
\begin{align}
    \int_\be\tau_h\bn\cdot\jump{u_h}~\ds &= \frac{1}{h}\|\jump{u_h}\|_{0,\be}^2~\forall~\be\in E_h, \qquad\int_\be\tau_h\bn\cdot u_h~\ds = \frac{1}{h}\|u_h\|_{0,\be}^2~\forall~\be\subset\Gamma_D, \\
                                         &\int_K \tau_h:\varepsilon(u_h)~\dx = -\|\varepsilon(u_h)\|_{0,K}^2~\forall~K\in\mathcal{T}_h.
\end{align}
Using~\eqref{eq:b_expanded},
\begin{equation}
    b(\tau_h,u_h) = -\sum_K\left(-\|\varepsilon(u_h)\|_{0,K}^2\right) + \sum_{\be\in E_h}\frac{1}{h}\|\jump{u_h}\|_{0,\be}^2 + \sum_{\be\subset\Gamma_D}\frac{1}{h}\|u_h\|_{0,\be}^2 = \|u_h\|_h^2.
\vspace{-1mm}\end{equation}
It remains to show that $\|u_h\|_h\gtrsim\verth{\tau_h}$. For every $\be\subset\Gamma_N$, since the DOFs associated to an edge and its endpoints determine $\tau_h\bn$ on that edge, we have $\tau_h\bn = 0$ by~\eqref{eq:tau_vertex},~\eqref{eq:tau_GammaN}, so it suffices to show $\|u_h\|_h\gtrsim\|\tau_h\|_{0,\Omega}$. 
By~\eqref{eq:tau_Eh}, 
~\eqref{eq:tau_internal}, we have $\pi_\be(\tau_h\bn) = \frac{1}{h}\jump{u_h}, \pi_K(\tau_h) = -\varepsilon(u_h)~\forall~\be\in E_h, \be\subset\Gamma_D, K\in\mathcal{T}_h$, where $\pi_\be:L^2(\be;\mathbb{R}^2)\to\mathcal{P}_1(\be;\mathbb{R}^2),\pi_K:L^2(K;\mathbb{S})\to\mathbb{S}$ are orthogonal projections. By equivalence of norms on $\widehat{\Sigma}_h$ on the reference cell $\hat{K}$, we have $\|\htau_h\|_{0,\hat{K}}^2 \lesssim \sum_{\text{vertices }\bhx\in\hat{K}}|\htau_h(\bhx)|^2 + \sum_{\bhe\subset\partial\hat{K}}\|\pi_{\bhe}(\htau\bhn)\|_{0,\bhe}^2 + \|\pi_{\hat{K}}\htau\|_{0,\hat{K}}^2$, so by a scaling argument we obtain $\|\tau_h\|_{0,K}^2\lesssim \sum_{\be\subset\partial K}\frac{1}{h}\|\jump{u_h}\|_{0,\be}^2 + \|\varepsilon(u_h)\|_{0,K}^2$.
This shows
\vspace{-1mm}\begin{equation}
    \inf_{u_h\in V_h}\sup_{\tau_h\in\Sigma_h}\frac{b(\tau_h,u_h)}{\verth{\tau_h}\|u_h\|_h}\gtrsim 1.
\vspace{-1mm}\end{equation}
\vspace{-3mm}\end{proof}

\vspace{-2mm}\begin{remark}\label{rem:normalproj}
    In analogy with the RT or BDM spaces considered in~\cite{Burman2021, Stenberg2011}, we have the useful equilibrium property $\Div\Sigma_h \subset V_h$, but in disanalogy we in general have $(\Sigma_h\vert_{\be})\bn\not\subset V_h\vert_\be$ on edges $\be$.
\end{remark}

For error estimation, we employ
intermediate approximants 
$\Pi_h\sigma, P_hu$, 
where $P_h:L^2(\Omega;\mathbb{R}^2)\to V_h$ denotes the orthogonal projection and $\Pi_h:H^1(\Omega;\mathbb{S})\to\Sigma_h$ denotes the 
Cl\'ement-like 
interpolation operator constructed by Arnold and Winther~\cite{arnold2002}, which enjoy the following approximation properties for all $(\tau, v)\in H^1(\Omega;\mathbb{S})\times L^2(\Omega;\mathbb{R}^2)$:
\begin{enumerate*}[label=(\roman*)]
    \item\label{item:comm} $\Div\Pi_h\tau = P_h\Div\tau$ (as in a smoothed variant of the commuting diagram~\eqref{eq:AW_complex});
    \item\label{item:edge_comm} $(\Pi_h\tau)\bn = \pi_\be(\tau\bn)$ on all edges $\be$;
    \item\label{item:Pi} $\|\tau - \Pi_h\tau\|_{0,\Omega}\lesssim h^m\|\tau\|_{m,\Omega}, 1\leq m\leq 3$;
    \item\label{item:P} $\|v - P_hv\|_{0,\Omega}\lesssim h^m\|v\|_{m,\Omega}, 0\leq m\leq 2$.
\end{enumerate*}
\vspace{-3mm}\begin{lemma}\label{lem:H2_approx}
    For each $K\in\mathcal{T}_h$ and $\tau\in H^2(\Omega;\mathbb{S})$, we have
    \begin{equation}
        \|\tau - \Pi_h\tau\|_{1,K}\lesssim h\|\tau\|_{2,S_K},
    \end{equation}
    where $S_K$ is a patch of cells neighbouring $K$ such that $\{S_K\}_{K\in\mathcal{T}_h}$ has the finite overlapping property.
\vspace{-2mm}\end{lemma}
\begin{proof}
    We adapt the proof 
    of~\ref{item:Pi} 
    in~\cite{arnold2002}, from which we recall that the error in $\Pi_h$ may be written as $I - \Pi_h = (I - \Pi^0_h)(I - R_h)$, where $R_h:L^2(\Omega;\mathbb{S})\to\Sigma_h\cap H^1(\Omega;\mathbb{S})$ is a Cl\'ement interpolant satisfying
    \begin{equation}\label{eq:Clement}
        \|\tau - R_h\tau\|_{j,K}\lesssim h^{m-j}\|\tau\|_{m,S_K}, \qquad 0\leq j\leq 1, \quad j\leq m\leq 3,
    \end{equation}
    and $S_K$ is a patch of the required form,
    and $\Pi^0_h:H^1(\Omega;\mathbb{S})\to\Sigma_h$ is the canonical interpolation operator except at the vertices, at which $(\Pi^0_h\tau)(\bx) \coloneqq 0~\forall~\bx\in\Lambda_h$. It is easily checked that $\Pi^0_{\hat{K}}$, the restriction of $\Pi^0_h$ to a single cell $\hat{K}$, is bounded from $H^1(\hat{K};\mathbb{S})$ to $H^r(\hat{K};\mathbb{S})$ for all $r\geq 0$; choosing $r = 1$, by a scaling argument we obtain
    $
        \|\Pi^0_K\tau\|_{1,K}\lesssim h^{-1}\|\tau\|_{0,K} + \|\tau\|_{1,K},
    $
    so
    \begin{equation}
        \|\Pi^0_h(I - R_h)\tau\|_{1,K}\lesssim h^{-1}\|(I - R_h)\tau\|_{0,K} + \|(I - R_h)\tau\|_{1,K}\lesssim h\|\tau\|_{2,S_K},
    \end{equation}
    which gives the result when combined with~\eqref{eq:Clement}.
\vspace{-3mm}\end{proof}

We make an additional regularity assumption on the stress field for error analysis.

\vspace{-1mm}\begin{proposition}\label{prop:nitsche_error}
    Let $(\sigma,u)\in H^2(\Omega;\mathbb{S})\times H^2(\Omega;\mathbb{R}^2)$.
    We have the error estimate
    \vspace{-1mm}\begin{equation}
        \|\sigma - \sigma_h\|_0 \lesssim h(\gamma\|\sigma\|_{2,\Omega} + \|u\|_{2,\Omega}).
    \vspace{-1mm}\end{equation}
\end{proposition}
\noindent In particular, we expect the traction residual to converge as $\|\sigma_h\bn - g\|_{0,\Gamma_N} = \mathcal{O}(h^{\frac{1}{2}})$.

\begin{proof}
By the Babu\v ska condition associated with the well-posed system~\eqref{eq:H-R_nitsche}, applied to $(\sigma_h - \Pi_h\sigma, u_h - P_hu)\in\Sigma_h\times V_h$, there are $(\tau_h, v_h)\in\Sigma_h\times V_h$ with $\verth{\tau_h} + \|v_h\|_h\leq 1$ and
\begin{equation}
    \verth{\sigma_h - \Pi_h\sigma} + \|u_h - P_hu\|_h \lesssim a_{h,\gamma}(\sigma_h - \Pi_h\sigma, \tau_h) + b(\tau_h, u_h - P_hu) + b(\sigma_h - \Pi_h\sigma, v_h)
\end{equation}
which by consistency of the system~\eqref{eq:H-R_nitsche} is equal to
\vspace{-1mm}\begin{align}
    a_{h,\gamma}&(\sigma - \Pi_h\sigma,\tau_h) + b(\tau_h, u - P_hu) + b(\sigma - \Pi_h\sigma,v_h)\\
                &= a_{h,\gamma}(\sigma - \Pi_h\sigma,\tau_h) - \int_{\Gamma_N}\tau_h\bn\cdot(u - P_hu)\ds\qquad\qquad\text{ using~\ref{item:comm},~\ref{item:edge_comm}}\\
                &= \int_\Omega \mathcal{A}(\sigma - \Pi_h\sigma):\tau_h~\dx + \underbrace{\frac{\gamma}{h}\int_{\Gamma_N}(\sigma - \Pi_h\sigma)\bn\cdot\tau_h\bn~\ds}_{(\ast)} - \int_{\Gamma_N}\tau_h\bn\cdot(u - P_hu)\ds = (\dagger).
\vspace{-1mm}\end{align}
By~\ref{item:edge_comm}, the term $(\ast)$ would vanish were it not for Remark~\ref{rem:normalproj}. Employing multiplicative trace inequalities, approximation properties of $P_h$, and Lemma~\ref{lem:H2_approx}, we have
\vspace{-1mm}\begin{align}
    (\dagger) &\lesssim \|\sigma - \Pi_h\sigma\|_0\|\tau_h\|_0 + \|\tau_h\|_{0,\Gamma_N}\left(\frac{\gamma}{h}\|(\sigma - \Pi_h\sigma)\bn\|_{0,\Gamma_N} + \|u - P_hu\|_{0,\Gamma_N}\right)\\
              &\lesssim h^2\|\sigma\|_{2,\Omega} + h^{\frac{1}{2}}\sum_{\be\subset\Gamma_N}\left(\frac{\gamma}{h}\|\sigma - \Pi_h\sigma\|_{0,\be} + \|u - P_hu\|_{0,\be}\right)\\
              &\lesssim h^2\|\sigma\|_{2,\Omega} + h^{\frac{1}{2}}\sum_K\left(\frac{\gamma}{h}\|\sigma - \Pi_h\sigma\|_{0,K}^{\frac{1}{2}}\|\sigma - \Pi_h\sigma\|_{1,K}^{\frac{1}{2}} + \|u - P_hu\|_{0,K}^{\frac{1}{2}}\|u - P_hu\|_{1,K}^{\frac{1}{2}}\right)\\
              &\lesssim h^2\|\sigma\|_{2,\Omega} + h^{\frac{1}{2}}\left( \frac{\gamma}{h}h\|\sigma\|_{1,\Omega}^{\frac{1}{2}}h^{\frac{1}{2}}\|\sigma\|_{2,\Omega} + h\|u\|_{2,\Omega}^{\frac{1}{2}}h^{\frac{1}{2}}\|u\|_{2,\Omega}^{\frac{1}{2}}\right) = h^2\|\sigma\|_{2,\Omega} + \gamma h\|\sigma\|_{2,\Omega} + h^2\|u\|_{2,\Omega},
\vspace{-1mm}\end{align}
so
\vspace{-1mm}\begin{align}
    \|\sigma - \sigma_h\|_0 &\leq \|\sigma - \Pi_h\sigma\|_0 + \|\sigma_h - \Pi_h\sigma\|_0\lesssim h^2\|\sigma\|_2 + \verth{\sigma_h - \Pi_h\sigma} + \|u_h - P_hu\|_h\\
                            &\lesssim h^2\|\sigma\|_2 + \gamma h\|\sigma\|_{2}+ h^2\|u\|_{2}.
\vspace{-1mm}\end{align}
\vspace{-3mm}\end{proof}

\vspace{-3mm}\section{An exterior calculus perspective}\label{sec:FEEC}
\vspace{-2mm}

The finite element exterior calculus~\cite{AFW2006} is a framework for structure-preservation in finite element discretizations, and for constructing finite element spaces as subcomplexes of complexes of function spaces of differential forms. The inventions of the MTW and AW elements we have considered were motivated by the Stokes~\eqref{eq:stokes_complex} and elasticity~\eqref{eq:AW_complex} complexes respectively. In this Section, we consider the application of FEEC to Piola transformation theory, and to the construction of multigrid smoothers.

\vspace{-3mm}\subsection{Uniform construction of the pullbacks}\label{sec:pullbacks}
\vspace{-2mm}

The Piola transforms~\eqref{eq:single_Piola}--\eqref{eq:double_cov_Piola} may be regarded as the analogy of the standard pullback~\eqref{eq:pullback} for $\hdiv$- and $\hcurl$-based spaces, but in fact the pullbacks may be defined uniformly in a manner guided by the FEEC~\citep[p.~35]{lithesis},~\citep[Section 6.2.5]{Arnold_FEEC}; we here employ terminology for which we refer the reader to~\citep[Ch.~6]{Arnold_FEEC}.\footnote{We remark that much of the literature employing Piola transforms does little to motivate their definition, beyond stating the \textit{consequences} of the definitions.}

Let us regard the physical and reference cells $K,\hat{K}$ as 
submanifolds of dimension $d$ in $\mathbb{R}^d$, 
with $F:\hat{K}\to K$ a diffeomorphism. 
Denote by $\alt^k\mathbb{R}^d$ the space of alternating $k$-linear forms on $\mathbb{R}^d$, and
for 
$M\in\{K,\hat{K}\}$ let $\Lambda^k(M)$ denote the space of differential $k$-forms on $M$, of which functions in the Sobolev spaces we consider will be scalar, vector, and tensor proxies, and the operators of vector calculus will be proxies for the exterior derivative $\mathrm{d}:\Lambda^k\to\Lambda^{k+1}$. Scalar fields may be identified with $0$-forms or $d$-forms, and vector fields with $1$-forms or $(d-1)$-forms. 
We may specify $L^2$-integrability of differential form coefficients with $L^2\Lambda^k(\Omega) \coloneqq L^2(\Omega;\alt^k\mathbb{R}^d)$; 
Sobolev spaces of differential forms may be defined as $H\Lambda^k \coloneqq \{\omega\in L^2\Lambda^k~\vert~\mathrm{d}\omega\in L^2\Lambda^{k+1}\}$, and correspond to the conventional spaces of vector calculus via $H\Lambda^0\simeq H^1, H\Lambda^1\simeq\hcurl, H\Lambda^{d-1}\simeq\hdiv$, and $H\Lambda^d\simeq L^2$.\footnote{We already employ $\mathbb{R}^d$ as a proxy for the tangent spaces $\{T_\bx M\}_{\bx\in M}$ on which the alternating forms are defined.}

An elastic stress field $\mathcal{T}$ on $M$ may naturally be identified with 
an $(\alt^{d-1}\mathbb{R}^d)$-valued $(d-1)$-form $\tau\in\Lambda^{d-1}(M;\alt^{d-1}\mathbb{R}^d)$, which by the proxy of $\mathbb{R}^d$ for $\alt^{d-1}\mathbb{R}^d$ may be identified with $\Lambda^{d-1}(M;\mathbb{R}^d)$,
since when integrated over a codimension-1 submanifold (such as the boundary of a subdomain), it should give a vector representing force~\cite{AFW2006}\citep[p.~618]{Frankel}. Applying the Hodge star gives an element $\star\tau\in\Lambda^1(M;\mathbb{R}^d)$, i.e.~a linear map $\mathbb{R}^d\to\mathbb{R}^d$ (hence, a matrix) at every point of $M$, which is the classical characterization of stress. 
Alternatively, $\mathcal{T}$ may 
be identified with a (symmetric) \textit{co}variant 2-tensor field in $\Lambda^1(M;\alt^1\mathbb{R}^d)$, which to each point assigns a (symmetric) bilinear form on $\mathbb{R}^d$~\citep[p.~10]{lithesis}.

Given $\hat{\omega}\in\Lambda^k(\hat{K})$, the derivative of the inverse diffeomorphism $J^{-1}(\bx) = (F^{-1})'(\bx)$ at $\bx\in K$ induces an element $\omega = J^{-*}\hat{\omega}\in\Lambda^k(K)$ pointwise via $(J^{-*}\hat{\omega})_\bx \coloneqq J^{-1}(\bx)^*\hat{\omega}_{F^{-1}(\bx)}$, the \textit{pullback} of $\hat{\omega}$ under $F^{-1}$, where$~^*$ denotes the algebraic pullback of a linear map $L:\mathbb{R}^d\to\mathbb{R}^d$ given by $L^*\eta(v_1,\ldots,v_k) \coloneqq \eta(Lv_1,\ldots,Lv_k)$ for $\eta\in\Lambda^K(M), v_i\in\mathbb{R}^d$. For scalar fields representing $0$-forms $\hat{\omega}\in H\Lambda^0(\hat{K})$ (i.e.~a constant map at each point), the scalar proxy for the resulting $0$-form $\omega$ is easily seen to be given simply by precomposition with $F^{-1}$, because $(J^{-*}\hat{\omega})_\bx = \hat{\omega}_{F^{-1}(\bx)}$, which for the proxy gives exactly the standard pullback~\eqref{eq:pullback}.
For a vector field $\hat{w}$ representing a $1$-form $\hat{\omega}\in H\Lambda^1(\hat{K})$, the canonical identification is $\hat{\omega} = \sum_{i=1}^d\hat{w}_i\mathrm{d}\bhx^i$, so that $\hat{\omega}_{\bhx}(v) = \hat{w}(\bhx)\cdot v$ (i.e.~$\hat{w}$ is the pointwise Riesz representative of $\hat{\omega}$), so
\vspace{-1mm}\begin{equation}
    (J^{-*}\hat{\omega})_\bx(v) = \hat{\omega}_{F^{-1}(\bx)}(J^{-1}(\bx)v) = \hat{w}(F^{-1}(\bx))\cdot(J^{-1}(\bx)v) = (J^{-T}(\bx)\hat{w}(F^{-1}(\bx)))\cdot v,
\vspace{-1mm}\end{equation}
hence the transformed proxy is given by $\mathcal{F}^{\curl}(\hat{w})$. A vector field $\hat{w}$ can alternatively represent a $(d-1)$-form via
\vspace{-1mm}\begin{align}
    \hat{\omega}_{\bhx}(v_1,\ldots,v_{d-1}) = \det[\hat{w}(\bhx)\vert v_1\vert \ldots\vert v_{d-1}] =: \sum_{i=1}^d \hat{w}_i(\bhx)\phi_i(v_1,\ldots,v_{d-1}) = \hat{w}(\bhx)\cdot\phi(v_1,\ldots,v_{d-1}),\label{eq:determinant_minors}
\vspace{-1mm}\end{align}
for which
\vspace{-1mm}\begin{equation}
\begin{aligned}
    (J^{-*}\hat{\omega})_\bx(v_1,\ldots,v_{d-1}) = \det[\hat{w}(F^{-1}(\bx))\vert J^{-1}v_1\vert \ldots\vert J^{-1}v_{d-1}] &= \det(J^{-1})\det[J\hat{w}(F^{-1}(\bx))\vert v_1\vert \ldots\vert v_{d-1}] \\
    &= \det[\Fdiv(\hat{w})(\bx)\vert v_1\vert \ldots\vert v_{d-1}].
\end{aligned}
\vspace{-1mm}\end{equation}
The derivation of the tensor transforms are analogous; a (symmetric) covariant 2-tensor field $\htau\in\Lambda^1(\hat{K};\alt^1\mathbb{R}^d)$ has as proxy a (symmetric) matrix field $\hat{\mathcal{T}}$ via the identification $\htau_{\bhx}(v_1,v_2) = v_1^T\hat{\mathcal{T}}(\bhx)v_2$, so that $(J^{-*}\htau)_\bx(v_1,v_2) = v_1^TJ^{-T}\hat{\mathcal{T}}(F^{-1}(\bx))J^{-1}v_2 = v_1^T\mathcal{F}^{\curl,\curl}(\hat{\mathcal{T}})(\bx)v_2$.
A (symmetric) matrix field $\hat{\mathcal{T}}$ may instead serve as a proxy for a (symmetric) 
form $\htau\in\Lambda^{d-1}(\hat{K};\alt^{d-1}\mathbb{R}^d)$ via the identification
\vspace{-1mm}\begin{align}
    \htau_{\bhx}(v_1,\ldots,v_{d-1})(z_1,\ldots,z_{d-1}) &= \phi(v_1,\ldots,v_{d-1})^T\hat{\mathcal{T}}(\bhx)\phi(z_1,\ldots,z_{d-1}), 
\vspace{-1mm}\end{align}
where $\phi$ is defined by~\eqref{eq:determinant_minors}. By definition of $\phi$ it is easily checked that $\phi(J^{-1}v_1,\ldots,J^{-1}v_{d-1}) = \frac{1}{\det J}J^T\phi(v_1,\ldots,v_{d-1})$, and hence
\vspace{-1mm}\begin{equation}
    \begin{aligned}
        (J^{-*}\htau)_{\bx}(v_1,\ldots,v_{d-1})(z_1,\ldots,z_{d-1}) &= \phi(J^{-1}v_1,\ldots,J^{-1}v_{d-1})^T\hat{\mathcal{T}}(F^{-1}(\bx))\phi(J^{-1}z_1,\ldots,z_{d-1}) \\
                                                 &= \phi(v_1,\ldots,v_{d-1})^T\Fdivdiv(\hat{\mathcal{T}})(\bx)\phi(z_1,\ldots,z_{d-1}).
    \end{aligned}
\vspace{-1mm}\end{equation}
Note that neither tensor transform is given by the row-wise application of the corresponding vector transform.

We hereafter assume that $J$ is spatially constant, so that the pullback commutes with the exterior derivative. 
In this case, a crucial property is that moments of the exterior derivatives against appropriate fields are 
preserved, a primary motivation for the use of the Piola pullbacks in mixed finite elements~\citep[Lemmas 2.1.6, 2.1.9]{BBF}. Let $\omega\in H\Lambda^k(K), \mu\in H\Lambda^{d-k-1}(K)$ with pullbacks $\hat{\omega},\hat{\mu}$, then since the pullback respects the exterior product and commutes with the exterior derivative,
\vspace{-1mm}\begin{equation}\label{eq:preserve_moments}
    \int_K\mu\wedge\mathrm{d}\omega = \int_K J^{-*}\hat{\mu}\wedge\mathrm{d}J^{-*}\hat{\omega} = \int_K J^{-*}\hat{\mu}\wedge J^{-*}\mathrm{d}\hat{\omega} = \int_K J^{-*}(\hat{\mu}\wedge\mathrm{d}\hat{\omega}) = \pm\int_{\hat{K}}\hat{\mu}\wedge\mathrm{d}\hat{\omega}.
\vspace{-1mm}\end{equation}
\noindent Note that since we have chosen to allow $F$ to reverse orientation, in general the moment preservation~\eqref{eq:preserve_moments} will only be up to sign. 
By Stokes' theorem, it also follows that
\vspace{-1mm}\begin{equation}\label{eq:preserve_boundary_moments}
    \int_{\partial K} \mu\wedge\omega = \pm\int_{\partial\hat{K}}\hat{\mu}\wedge\hat{\omega},
\vspace{-1mm}\end{equation}
\noindent where the boundary integrands are meant in the trace sense. 

Combining these identities with the characterization of constrained finite element spaces by kernels of constraint functionals, as in subsection~\eqref{sec:constraints}, also allows for a more elegant verification of their preservation by the Piola maps;
the constraints for the MTW and AW\textsuperscript{c} elements are preserved in light of~\eqref{eq:preserve_moments}, while the preservation of boundary constraints follows from the identity~\eqref{eq:preserve_boundary_moments} 
applied to appropriate $L^2(\partial K)$-orthogonal polynomials, as constructed e.g.~in~\citep[Section 3]{Olver2019}.

That the pullback commutes with the exterior derivative 
also implies the preservation of the kernel of the operators defining the spaces, a property implicit in the classical presentation of the Piola operators~(e.g.~\citep[Section 3.9]{Monk}). For example, if $K,\hat{K}\subset \mathbb{R}^2$ (respectively, $\mathbb{R}^3$), note that if $\hat{\phi}\in H^1(\hat{K})$ and $\phi \in H^1(K)$ is its standard pullback given by~\eqref{eq:pullback}, then by the chain rule 
we have $\nabla \phi = J^{-T}\hat{\nabla}\hat{\phi}$.
Since $\rot\nabla\equiv 0$ (resp. $\curl\nabla\equiv 0$), certainly both $\hat{\nabla}\hat{\phi}\in H(\rot,\hat{K})$ (resp. $H(\curl,\hat{K})$) and $\nabla \phi\in H(\rot,K)$ (resp. $H(\curl,K)$), and the covariant map~\eqref{eq:covariant} ``transform[s] vectors of $H(\curl;\Omega)$ like gradients''~\citep[p.~61]{BBF}, which, on simply connected domains, form the kernel of $\rot$ or $\curl$ by exactness of an appropriate de Rham sequence; hence
\vspace{-1mm}\begin{equation}
    \ker(\rot) = \mathcal{F}^{\curl}\left(\ker\left(\widehat{\rot}\right)\right)\quad\text{ in 2D }\qquad\text{ and }\qquad\ker(\curl) = \mathcal{F}^{\curl}\left(\ker\left(\widehat{\curl}\right)\right)\quad\text{ in 3D.}
\vspace{-1mm}\end{equation}
Analogously, by a direct computation,
\vspace{-1mm}\begin{equation}\label{eq:curl_Piola}
    \curl \phi = \frac{1}{\det J}J\widehat{\curl}~\hat{\phi},
\vspace{-1mm}\end{equation}
which is reflected in~\eqref{eq:single_Piola}. We now extend the discussion in~\citep[Section 3.9]{Monk} to the tensor-valued case in 2D. Let us define the matrix-valued $\curl$ of a vector field in 2D,
\vspace{-1mm}\begin{equation}
    \text{\textbf{curl}}~v \coloneqq \begin{bmatrix} \dfrac{\partial v_1}{\partial y} & -\dfrac{\partial v_1}{\partial x} \\ \dfrac{\partial v_2}{\partial y} & -\dfrac{\partial v_2}{\partial x} \end{bmatrix},
\vspace{-1mm}\end{equation}
for which
\vspace{-1mm}\begin{equation}\label{eq:mat_curl_Piola}
    \text{\textbf{curl}}~v = \frac{1}{\det J}\widehat{\text{\textbf{curl}}}~\hat{v}~J^T,
\vspace{-1mm}\end{equation}
and the Airy stress function associated with a scalar potential $\phi$
\vspace{-1mm}\begin{equation}\label{eq:airy_stress}
    \airy\phi \coloneqq \text{\textbf{curl}}\curl\phi = \begin{bmatrix} \dfrac{\partial^2\phi}{\partial y^2} & -\dfrac{\partial^2\phi}{\partial x\partial y} \\ -\dfrac{\partial^2\phi}{\partial x\partial y} & \dfrac{\partial^2\phi}{\partial x^2} \end{bmatrix},
\vspace{-1mm}\end{equation}
which is identically symmetric and divergence-free. Then for $\hat{\phi}\in H^2(\hat{K})$ with pullback $\phi\in H^2(K)$,
\vspace{-1mm}\begin{equation}
    \airy\phi = \frac{1}{(\det J)^2}J\Big(\widehat{\airy}~\hat{\phi}\Big)J^T = \Fdivdiv\left(\widehat{\airy}\hat{\phi}\right),
\vspace{-1mm}\end{equation}
so by exactness of the 2D stress complex~\eqref{eq:AW_complex}, we obtain that $\ker(\Div) = \Fdivdiv(\ker(\widehat{\Div}))$.
Similarly, for the $\rot$ of a symmetric matrix field applied row-wise, we calculate that $\nabla^2\phi = \mathcal{F}^{\curl,\curl}(\hat{\nabla}^2\hat{\phi})$, where $\nabla^2 H^2$ forms the kernel of $\rot$ on $H(\rot;\mathbb{S})$ 
by exactness of the 2D Hessian complex (e.g.~\citep[Remark 3.16]{pauly2020}).

This observation of kernel preservation may also be connected to the (albeit trivial) topologies of $K$ and $\hat{K}$: the pullbacks are isomorphisms between appropriate complexes on $K$ and $\hat{K}$, which moreover commute with $\mathrm{d}$, hence are cochain maps which preserve the cohomology of each domain.

\vspace{-3mm}\subsection{Kernel-capturing and robust multigrid}\label{sec:kernel_cap}
\vspace{-2mm}

It is now well-established that the characterization of the kernels of discretized differential operators is crucial for the design of robust multigrid schemes~\cite{FMSW_N-S, SchoeberlThesis, AFW2000, FKMW}; for both the MTW and AW-type elements, this is given precisely by their positions in the discrete exact complexes~\eqref{eq:MTW_complex} and~\eqref{eq:AW_complex}, as we now explain.
We describe multigrid relaxation in the framework of subspace correction methods~\cite{Xu1992}. Given a finite-dimensional Hilbert space $V$ of functions on a mesh, consider a decomposition
\vspace{-1mm}\begin{equation}\label{eq:schwarz_decomp}
    V = \sum_i V_i,
\vspace{-1mm}\end{equation}
where the sum need not be direct. The variational problem to be solved over $V$ often takes the form of a symmetric, coercive operator, perturbed by a positive semidefinite singular operator (such as a discretized divergence) which is scaled by some parameter $\alpha > 0$, a physical or penalty parameter which, as it increases, renders the problem difficult to solve. The seminal work of Sch\"oberl~\citep[Theorem 4.1]{SchoeberlThesis} revealed sufficient conditions for $\alpha$-robustness of the parallel subspace correction preconditioner induced by the decomposition~\eqref{eq:schwarz_decomp}; a key insight is that, if $\mathcal{N}$ denotes the kernel of the singular operator, the subspaces should satisfy the \textit{kernel-capturing property}
\vspace{-1mm}\begin{equation}\label{eq:kernel_decomp}
    \mathcal{N} = \sum_i(V_i\cap\mathcal{N}).
\vspace{-1mm}\end{equation}
For an $\hdiv$-based space $V$ such as the MTW velocity space in the Stokes complex~\eqref{eq:stokes_complex} or the AW stress space in the elasticity complex~\eqref{eq:AW_complex}, one choice of relaxation method is given by the \textit{vertex star iteration}~\cite{FMW}, which is induced by the subspaces
\vspace{-1mm}\begin{equation}\label{eq:vertex_star}
    V_i = \{v\in V~|~\supp(v)\subset K_i\},
\vspace{-1mm}\end{equation}
where $K_i$ denotes the patch of cells in the mesh sharing vertex $i$. Given $v\in V$ with $\Div v = 0$, where now $V\in\{V_h,\Sigma_h\}$ denotes one of the canonical MTW or AW spaces respectively, by exactness of the associated discrete complex when $\Omega$ is simply connected, we may write $\mathcal{C}\phi = v$ for some potential $\phi\in W$, where $(W, \mathcal{C})$ denote $(W_h/\mathbb{R}, \curl)$ or $(Q_h, \airy_h)$ respectively. Let now $\{\phi_i\}_i$ denote a basis for $W$, and write $\phi = \sum_i c_i\phi_i$. Then a divergence-free decomposition of $v$ is given by $v_i = c_i\mathcal{C}\phi_i$, since $v = \sum_i v_i$ and $v_i\in\mathcal{N}$ for each $i$, so it suffices to find subspaces $V_i$ with $\mathcal{C}\phi_i\in V_i$ for each $i$. That the vertex star~\eqref{eq:vertex_star} fulfils this property follows from 
inspection of the basis functions of $W$.

In general, Sch\"oberl's hypotheses also require that the splitting~\eqref{eq:schwarz_decomp} be stable in the $V$-norm and that the kernel splitting~\eqref{eq:kernel_decomp} be stable in the energy norm induced by the Galerkin projection of the coercive form, which does not follow automatically from the exactness of the discrete complex; typically, such bounds hold for the infinite-dimensional spaces, and hold also for the discrete complex if bounded commuting projections exist.

\vspace{-3mm}\section{Smoothers}\label{sec:smoothers}
\vspace{-2mm}

To demonstrate the effectiveness of our mapping techniques in practical numerical simulations, it is necessary to consider preconditioners for the partial differential equations discretized by the MTW and AW elements, and to exhibit composability of our implementation with 
the associated software stack~\cite{kirby2018solver}. 
As indicated in the previous subsection, the application of patch-based multigrid smoothers is natural for the $\hdiv$-discretizing elements we consider in this paper.

Building on the work of Benzi and Olshanskii~\cite{Benzi} and Sch\"oberl~\cite{SchoeberlThesis} among others, Farrell and coauthors have successfully developed parameter- and mesh-robust 
preconditioners
of augmented Lagrangian (AL) type, with highly specialized multigrid algorithms, for a host of nonlinear PDEs with saddle point structure~\cite{FMW, FMSW_N-S, GF, Gazca-Orozco_anisothermal,XFW,RRB_mixed,FMSW_elast,Laakmann2021}. 
We illustrate the AL method for AW elements, the application to MTW being analogous.
It consists of augmenting the Hellinger--Reissner Lagrangian~\eqref{eq:H-R_fnl} with a penalty term
\vspace{-1mm}\begin{equation}\label{eq:AL_energy}
    \mathcal{H}_\alpha(\sigma, u) \coloneqq \mathcal{H}(\sigma, u) + \frac{\alpha}{2}\int_\Omega \|\Div\sigma - f\|^2\dx
\vspace{-1mm}\end{equation}
for $\alpha \geq 0$. The AL term penalizes deviation from the set constrained by the momentum balance~\eqref{eq:equil}, but does not change the exact solution. 
Its more significant benefit is that it allows the control of the Schur complement of the discretized system, as we now explain. The stationarity condition of the augmented energy~\eqref{eq:AL_energy} gives rise to a saddle point system for the stress-displacement pair, which, by a block factorization, admits the well-known solution formula~\citep[p.~14]{Benzi2005}
\vspace{-1mm}\begin{equation}
    \begin{pmatrix}
        A_\alpha & B^T \\
        B & \\
    \end{pmatrix}^{-1}
    =
    \begin{pmatrix}
        \mathbb{I} & -A_\alpha^{-1}B^T \\
                    & \mathbb{I} \\
    \end{pmatrix}
    \begin{pmatrix}
        A_\alpha^{-1} &  \\
               & S_\alpha^{-1}
    \end{pmatrix}
    \begin{pmatrix}
        \mathbb{I}         &  \\
        -BA_\alpha^{-1}  & \mathbb{I} \\
    \end{pmatrix},
\vspace{-1mm}\end{equation}
where $A, B^T, B$ are discrete 
compliance,
strain, 
and divergence operators respectively, $S = -BA^{-1}B^T$ is the (in general, dense) Schur complement, and the subscript $\alpha$ denotes the same quantities but of the augmented system (so that $A_0 = A$ etc.). We wish to precondition this system 
for FGMRES Krylov iterations. In analogy to the velocity-pressure 
Stokes 
problem, for which the Schur complement is spectrally equivalent to the viscosity-weighted pressure mass matrix~\cite{Wathen1994}, it can be shown (denoting by $M_u$ the displacement mass matrix) that
\vspace{-1mm}\begin{equation}\label{eq:schur}
    \tilde{S}^{-1} \coloneqq -\alpha M_u^{-1} \sim S_\alpha^{-1}
\vspace{-1mm}\end{equation}
serves as an approximate inverse to the augmented Schur complement, at least for fixed mesh size, with the approximation improving as $\alpha\rightarrow\infty$. Preconditioning the Schur complement with~\eqref{eq:schur} is however in tradeoff with the augmented stress solve $A_\alpha^{-1}$; the AL term in $A_\alpha$ has a large kernel consisting of the infinite-dimensional affine space of tensor fields with divergence $f$, rendering standard multigrid relaxation schemes ineffective. 
We propose the vertex star relaxation as an $\alpha$-robust multigrid algorithm for this block.

\vspace{-1mm}\begin{remark}
    All three of the finite element spaces we implement in this paper are non-nested under uniform refinement; 
    we employ the default prolongation operator of Firedrake, which involves (i) lossless projection of the coarse function onto a DG space of the same degree, (ii) lossless natural injection of the coarse projection, from the coarse grid to the fine DG space, (iii) projection from the fine DG space to the fine finite element space. 
    We conjecture that multigrid convergence behaviour observed 
    in Section~\ref{sec:numerics}
    may be improved by the construction of specialized prolongation operators for each element
    which preserve, at least approximately, the kernel of the discrete divergence, 
    in the style of~\cite{FMW,FMSW_elast} 
    (this being the other component of the parameter-robust multigrid framework of Sch\"oberl),
    but we do not pursue this here.
\vspace{-2mm}\end{remark}

\vspace{-2mm}\begin{remark}\label{rem:vertex_divfree}
    Let $\tau_h\in \text{AW\textsuperscript{c}}(K)$ denote a local AW\textsuperscript{c} basis function dual to a vertex DOF for a given cell $K$. Because all other nodes vanish at $\tau_h$, its full normal component vanishes along each edge, and its components have vanishing mean. Since $\varepsilon(\Div\tau_h)$ is a constant matrix,
    \vspace{-1mm}\begin{align}
        \|\Div\tau_h\|_{L^2(K)}^2 = \int_{\partial K} (\tau_h\bn)\cdot(\Div\tau_h)\ds - \int_K \tau_h:\varepsilon(\Div\tau_h)\dx = 0;
    \vspace{-1mm}\end{align}
    thus, compared to AW\textsuperscript{nc}, the nodal AW\textsuperscript{c} basis functions arising from the `extra' vertex DOFs contribute only to the divergence-free subspace. It follows that similar multigrid transfer operators may work well for both these elements.
\vspace{-2mm}\end{remark}

\vspace{-3mm}\section{Numerical examples}\label{sec:numerics}
\vspace{-2mm}

For the case of affine rather than Piola transformations, we have described the inclusion of this theory in the Firedrake code stack in~\cite{finat-zany}, and since this stack already understands Piola transformation the process is quite analogous.   We must implement each new reference element in FIAT~\cite{Kirby:2004} and wrap it into FInAT~\cite{Homolya2017finat}.  The FInAT wrapper also requires a function to construct abstract syntax for the basis transformation in terms of callbacks provided by the form compiler to obtain symbols for geometric quantities such as Jacobians, and physical and reference normal and tangent vectors.  The new elements must also be registered with UFL~\cite{Alnaes2014UFL} and \texttt{tsfc}~\cite{Homolya2018tsfc}.
We now consider several test problems to validate our implementation and demonstrate its capabilities.

Linear systems in the manufactured solution examples, and on the coarsest grid in the multigrid examples, were solved by sparse LU factorization with MUMPS~\cite{mumps} and PETSc~\cite{petsc}.

\vspace{-3mm}\subsection{Mardal--Tai--Winther}
\vspace{-2mm}

We consider the following parameterized saddle point system: given $f, g\in L^2(\Omega;\mathbb{R}^2), j\in H^{1/2}(\Gamma_D;\mathbb{R}^2)$, seek $(u, p)\in(H_{j\cdot\bn,\Gamma_D}(\Div)\cap\epsilon H^1_{j,\Gamma_D}(\Omega;\mathbb{R}^2))\times L^2(\Omega)$ such that
\vspace{-1mm}\begin{subequations}\label{eq:mtw_saddle}
    \begin{alignat}{2}
        \left( I - \epsilon^2 \Delta \right) u + \nabla p & = f &&\qquad\qquad\text{ in }\Omega,\\
        \Div u & = g &&\qquad\qquad\text{ in }\Omega,\\
        u &= j &&\qquad\qquad\text{ on }\Gamma_D,\label{eq:MTW_boundary_data}\\
        \epsilon^2\nabla u\bn - p\bn &= 0 &&\qquad\qquad\text{ on }\Gamma_N,
    \end{alignat}
\end{subequations}
in some domain $\Omega \subset \mathbb{R}^2$ with partitioned boundary $\Gamma_D\cup\Gamma_N$. 
Here, $\epsilon > 0$ is some positive parameter. Equations of this form arise in implicit timestepping for Stokes equations, for which $\epsilon$ corresponds to the square root of the timestep, and also in certain fractional-step methods for Navier--Stokes.  For moderate $\epsilon$, the equations are similar to the standard Stokes system and regular stable discretizations such as Taylor--Hood are readily applicable.  However, as $\epsilon \searrow 0$, the equations degenerate to an $\hdiv$ Darcy system (equivalently, a mixed Poisson problem) and the stability and error estimates of standard $H^1$ methods degrade.  The Mardal--Tai--Winther element, however, is stable for both $\hdiv$ and $H^1$ (though nonconforming) and for a scale of spaces in between; this is manifested in exhibiting stability and accuracy uniformly in $\epsilon$ for the system~\eqref{eq:mtw_saddle}.

\begin{figure}[H]
    \vspace{-3mm}
    \includegraphics[width=0.3\linewidth]{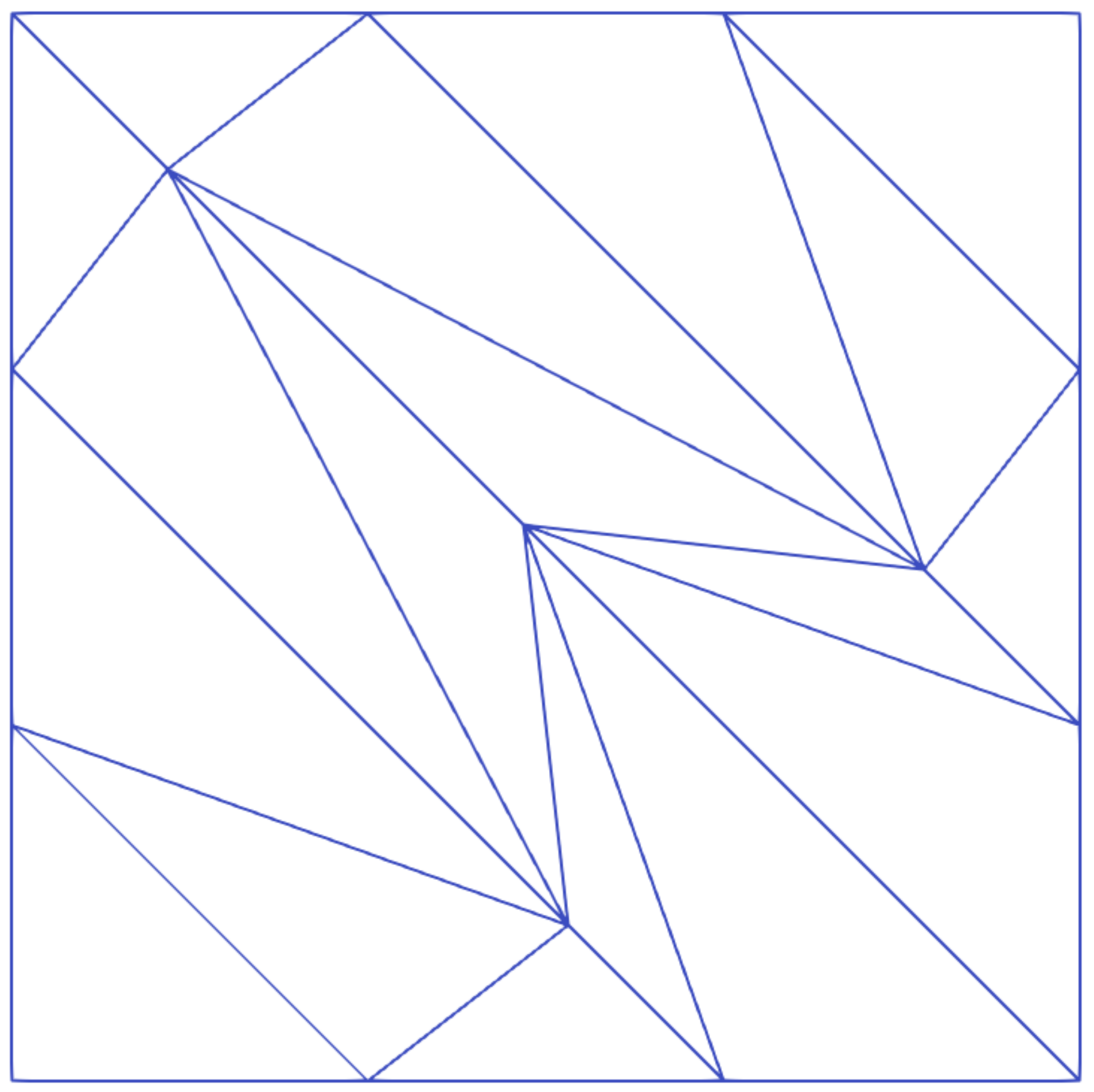}
    \vspace{-4mm}
    \caption{A warped mesh employed to check convergence under refinement.}\label{fig:warped}
    \vspace{-5mm}
\end{figure}

On the unit square $\Omega = (0,1)^2$, with $\Gamma_N = \{y=1\}, \Gamma_D = \partial\Omega\setminus\Gamma_N$, we approximate the manufactured solution $u(x,y) = (2^{1-y},0)^T,$ $p(x,y) = \cos(\pi x)\cos(2\pi y)$.
The essential boundary data~\eqref{eq:MTW_boundary_data} is enforced weakly using the boundary DOFs and an interpolant of the given boundary data $j$. Denoting the resulting global constrained MTW space by $V_{h,j}$, and the global DG$(0)$ pressure space by $Q_h$, we seek $(u_h, p_h)\in V_{h,j}\times Q_h$ such that
\vspace{-1mm}\begin{equation}\label{eq:MTW_MMS}
    \begin{aligned}
        \int_\Omega u_h\cdot v_h + \epsilon^2\nabla_h u_h:\nabla_h v_h - p_h\Div v_h~\dx &= \int_\Omega f\cdot v_h~\dx + \int_{\Gamma_N}K\bn\cdot v_h~\ds&&~\forall~v_h\in V_{h,0},\\
        -\int_\Omega q_h\Div u_h~\dx &= \int_\Omega g\cdot q_h~\dx &&~\forall~q_h\in Q_h,
    \end{aligned}
\vspace{-1mm}\end{equation}
where now $f,g,j,K$ are defined by appropriate residuals of the manufactured solution. In order to show that the mapping techniques apply on general unstructured meshes, we perturb 
the interior vertices of the coarsest mesh as pictured in Figure~\ref{fig:warped}; further refinements are obtained uniformly. In Tables~\ref{tab:MTW_MMS_velocity} and~\ref{tab:MTW_MMS_pressure}, we reproduce the approximately $\epsilon$-independent expected orders of convergence (EOCs) in both velocity and pressure proved in~\cite{mardal2002} for the MTW-DG$(0)$ pair, for a representative range of parameters $\epsilon$. Notice that the errors themselves are essentially $\epsilon$-independent, and \textit{a fortiori} so are the convergence rates, demonstrating the robustness of the element between these contrasting regimes. 

\begin{table}[H]
    \vspace{-3mm}
    \centering
    {\tiny
    \begin{tabular}{p{1.5cm}||p{1.5cm}|p{1.5cm}|p{1.5cm}|p{1.5cm}|p{1.5cm}|p{1.5cm}||p{1.5cm}}
                \toprule
                $\epsilon~\backslash~N$ &  $2^0$ & $2^1$ & $2^2$ & $2^3$ & $2^4$ & $2^5$ & EOC\\
                \midrule
                1          &  9.09$\times 10^{-3}$ & 2.45$\times 10^{-3}$ & 6.47$\times 10^{-4}$ & 1.68$\times 10^{-4}$ & 4.27 $\times 10^{-5}$ & 1.08 $\times 10^{-5}$ & 1.94\\
                $2^{-2}$       &  8.79$\times 10^{-3}$ & 2.39$\times 10^{-3}$ & 6.29$\times 10^{-4}$ & 1.62$\times 10^{-4}$ & 4.13 $\times 10^{-5}$ & 1.04 $\times 10^{-5}$ & 1.94\\
                $2^{-4}$   &  7.49$\times 10^{-3}$ & 2.17$\times 10^{-3}$  & 5.71$\times 10^{-4}$ & 1.46$\times 10^{-4}$ & 3.70 $\times 10^{-5}$ & 9.32 $\times 10^{-6}$  & 1.93\\
                $2^{-6}$   &  6.81$\times 10^{-3}$ & 1.93$\times 10^{-3}$ & 5.25$\times 10^{-4}$ & 1.36$\times 10^{-4}$ & 3.46 $\times 10^{-5}$ & 8.70 $\times 10^{-6}$ & 1.92\\
                $2^{-8}$ &  6.77$\times 10^{-3}$ & 1.88$\times 10^{-3}$ & 5.03$\times 10^{-4}$ & 1.31$\times 10^{-4}$ & 3.36 $\times 10^{-5}$ & 8.51 $\times 10^{-6}$ & 1.93\\
                $2^{-10}$ &  6.77$\times 10^{-3}$ & 1.88$\times 10^{-3}$ & 5.02$\times 10^{-4}$ & 1.30$\times 10^{-4}$ & 3.30 $\times 10^{-5}$ & 8.37 $\times 10^{-6}$ & 1.93\\
                0          &  6.77$\times 10^{-3}$ & 1.88$\times 10^{-3}$ & 5.02$\times 10^{-4}$ & 1.30$\times 10^{-4}$   & 3.30 $\times 10^{-5}$ & 8.32 $\times 10^{-6}$ & 1.93\\
                \bottomrule
            \end{tabular}
            }
    \caption{Errors and convergence rates in $L^2$ of the MTW velocity for the model problem~\eqref{eq:MTW_MMS}. Here and below, $N$ denotes the uniform refinement factor with respect to the original mesh.}\label{tab:MTW_MMS_velocity}
    \vspace{-2mm}
\end{table}


\begin{table}[H]
    \vspace{-3mm}
    \centering
    {\tiny
    \begin{tabular}{p{1.5cm}||p{1.5cm}|p{1.5cm}|p{1.5cm}|p{1.5cm}|p{1.5cm}|p{1.5cm}||p{1.5cm}}
            \toprule
            $\epsilon~\backslash~N$ &  $2^0$ & $2^1$ & $2^2$ & $2^3$ & $2^4$ & $2^5$ & EOC\\
            \midrule
            1          &  3.48$\times 10^{-1}$ & 1.83$\times 10^{-1}$ & 9.41$\times 10^{-2}$ & 4.73$\times 10^{-2}$ & 2.37$\times 10^{-2}$ & 1.18$\times 10^{-2}$ & 0.975\\
            $2^{-2}$       &  3.08$\times 10^{-1}$ & 1.79$\times 10^{-1}$  & 9.35$\times 10^{-2}$  & 4.71$\times 10^{-2}$ & 2.36$\times 10^{-2}$ & 1.18$\times 10^{-2}$ & 0.941\\
            $2^{-4}$     &  3.08$\times 10^{-1}$ & 1.79$\times 10^{-1}$ & 9.35$\times 10^{-2}$ & 4.71$\times 10^{-2}$ & 2.36$\times 10^{-2}$ & 1.18$\times 10^{-2}$  & 0.941\\
            $2^{-6}$   &  3.08$\times 10^{-1}$  & 1.79$\times 10^{-1}$ & 9.35$\times 10^{-2}$ & 4.71$\times 10^{-2}$ & 2.36$\times 10^{-2}$ & 1.18$\times 10^{-2}$  & 0.941\\
            $2^{-8}$ &  3.08$\times 10^{-1}$  & 1.79$\times 10^{-1}$ & 9.35$\times 10^{-2}$ & 4.71$\times 10^{-2}$ & 2.36$\times 10^{-2}$ & 1.18$\times 10^{-2}$  & 0.941\\
            $2^{-10}$ &  3.08$\times 10^{-1}$  & 1.79$\times 10^{-1}$ & 9.35$\times 10^{-2}$ & 4.71$\times 10^{-2}$ & 2.36$\times 10^{-2}$ & 1.18$\times 10^{-2}$  & 0.941\\
            0          &  3.08$\times 10^{-1}$  & 1.79$\times 10^{-1}$ & 9.35$\times 10^{-2}$ & 4.71$\times 10^{-2}$ & 2.36$\times 10^{-2}$ & 1.18$\times 10^{-2}$  & 0.941\\
            \bottomrule
        \end{tabular}
        }
    \caption{Errors and convergence rates in $L^2$ of pressure for~\eqref{eq:MTW_MMS}.}\label{tab:MTW_MMS_pressure}
    \vspace{-5mm}
\end{table}

\begin{figure}[H]
    \centering
    \includegraphics[width=0.6\linewidth]{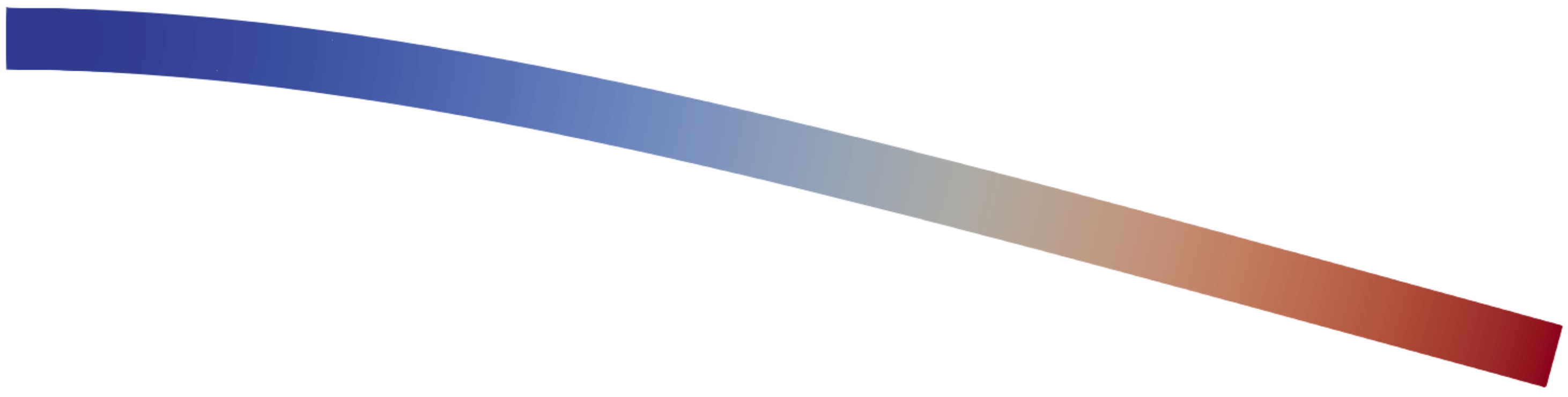}
    \vspace{-4mm}
    \caption{A nearly incompressible cantilever beam.}\label{fig:MTW_beam}
    \vspace{-5mm}
\end{figure}

We now demonstrate instead the robustness of our proposed multigrid solver for the MTW element applied to primal planar linear elasticity near the incompressible limit, a regime in which the MTW discretization will itself exhibit robust convergence; this example adapts the demonstration in~\cite{CometFEniCS}. 
We consider a simple cantilever beam $\Omega = [0,25]\times[0,1]$, fixed to a wall at the left-hand end $\Gamma_D$, stress-free on the rest of its boundary, in plane stress conditions, and subject only to the force of its own weight $f = (0,-10^{-3})^T$. Employing a Union-Jack-crossed mesh, and the associated MTW space $V_{h,0}$, we seek $u_h\in V_{h,0}$ such that
\vspace{-1mm}\begin{equation}\label{eq:MTW_primal_elast}
    \int_\Omega \mathbb{C}\varepsilon_h(u_h):\varepsilon_h(v_h)\dx = \int_\Omega f\cdot v_h~\dx~\forall~v_h\in V_{h,0},
\vspace{-1mm}\end{equation}
where the shear modulus of the material is fixed at $\mu = 3.8\times 10^4$ and the Poisson ratio is initialized at $\nu = 0.3$. In Table~\ref{tab:MTW_robust}, for a discretization of $7.23\times 10^5$ DOFs, we demonstrate robustness of the associated star relaxation, implemented via \texttt{PCPATCH}~\cite{FKMW}, as the Poisson ratio approaches the critical value $\frac{1}{2}$; the solution at highest Poisson ratio, enlarged for the sake of illustration, is pictured in Figure~\ref{fig:MTW_beam}. This was carried out 
with relative $\ell^2$-norm tolerance $10^{-5}$ for the outermost Krylov iterations.

\begin{table}[H]
    \centering
   {\tiny
       \begin{tabular}{l|l}
        \toprule
        Poisson ratio & Krylov iterations\\
        \midrule
         0.3 & 13\\
         0.45 & 14\\
         0.49 & 14\\
         0.4999 & 14\\
         0.4999999 & 14\\
        \bottomrule
        \end{tabular}
    }
    \caption{Moderate and approximately constant Krylov iteration counts for the solution of~\eqref{eq:MTW_primal_elast} near the incompressible limit, after preconditioning via the vertex star relaxation described in subsection~\ref{sec:kernel_cap}.}\label{tab:MTW_robust}
    \vspace{-5mm}
\end{table}


\vspace{-3mm}\subsection{The Hellinger--Reissner principle}
\vspace{-2mm}

We now consider the canonical Hellinger--Reissner problem~\eqref{eq:MMS_H-R} and verify the convergence results proved by Arnold and Winther~\citep[Theorem 5.1]{arnold2002}\citep[Theorem 4.1]{arnold2003}, summarized in Table~\ref{tab:convergence_rates}. The approximation order $m$ varies because in some cases, higher-order convergence may be obtained if the exact solution pair admits improved Sobolev regularity.

\begin{table}[H]
    \vspace{-2mm}
    \centering
    {\footnotesize
    \begin{tabular}{p{1.4cm}|p{1.9cm}|p{1.9cm}}
        \toprule
        Variable & AW\textsuperscript{c} & AW\textsuperscript{nc} \\
        \midrule
        $\sigma$     & $1\leq m\leq 3$ & 1               \\
        $\Div\sigma$ & $0\leq m\leq 2$ & $0\leq m\leq 2$               \\
        $u$          & $1\leq m\leq 2$ & 1               \\
        \bottomrule
    \end{tabular}
    }
    \caption{Ranges of approximation order $m$ by the AW elements of the elasticity variables in the $L^2$-norm.}\label{tab:convergence_rates}
    \vspace{-5mm}
\end{table}


On the unit square $\Omega = (0,1)^2$, we consider the 
displacement from Bramwell et al.~\cite{Bramwell}: 
    $u(x,y) = (\sin(\pi x)\sin(\pi y), \sin(\pi x)\sin(\pi y))^T,$
which satisfies an unphysical homogeneous boundary condition, 
and the stress field
\vspace{-1mm}\begin{equation}
    \sigma(x,y) = \begin{pmatrix} \cos(\pi x)\cos(3\pi y) & y + 2\cos(\frac{\pi}{2}x) \\ y + 2\cos(\frac{\pi}{2}x) & -\sin(3\pi x)\cos(2\pi x) \end{pmatrix}.
\vspace{-1mm}\end{equation}
In Tables~\ref{tab:first_AWnc_MMS} and~\ref{tab:last_AWnc_MMS}, we check convergence of AW\textsuperscript{nc} in the $L^2$-norms in the incompressible limit $\nu\nearrow\frac{1}{2}$, again on a perturbed mesh: 

\begin{table}[H]
    \vspace{-3mm}
    {\tiny
    \begin{tabular}{ p{0.5cm}|p{1.8cm}|p{1.8cm}|p{1.8cm}|p{1.5cm}|p{2.15cm}|p{2.2cm} }
        \toprule
        $N$ & $u$ error & $u$ EOC & $\sigma$ error & $\sigma$ EOC & $\Div_h(\sigma)$ error & $\Div_h(\sigma)$ EOC\\
        \midrule
          $2^0$ &  1.05$\times 10^{-1}$  &      --      &       5.51$\times 10^{-1}$  &      --       &          8.32$\times 10^{-1}$   &           --\\
        $2^1$ &  2.78$\times 10^{-2}$     &   1.91  &     2.54$\times 10^{-2}$    &    1.12   &        2.68$\times 10^{-1}$     &         1.64\\
        $2^2$ &  7.07$\times 10^{-3}$   &   1.98  &     1.18$\times 10^{-1}$    &    1.10   &        7.15$\times 10^{-2}$    &         1.91\\
        $2^3$ &  1.77$\times 10^{-3}$  &    2.00 &      5.76$\times 10^{-2}$   &     1.04  &         1.82$\times 10^{-2}$   &          1.98\\
        $2^4$ &  4.43$\times 10^{-4}$  &    2.00 &      2.85$\times 10^{-2}$  &     1.02  &         4.56$\times 10^{-3}$  &          1.99\\
        $2^5$ &  1.11$\times 10^{-4}$ &    2.00 &      1.42$\times 10^{-2}$  &     1.00  &         1.14$\times 10^{-3}$  &          2.00\\ 
        \bottomrule
    \end{tabular}
    }
    \caption{Errors and convergence rates in the $L^2$-norms using AW\textsuperscript{nc} for the model problem~\eqref{eq:MMS_H-R}, with $\nu = 0.25, \mu = 1$.}\label{tab:first_AWnc_MMS}
    \vspace{-5mm}
\end{table}


\begin{table}[H]
    \vspace{-3mm}
    {\tiny
    \begin{tabular}{ p{0.5cm}|p{1.8cm}|p{1.8cm}|p{1.8cm}|p{1.5cm}|p{2.15cm}|p{2.2cm} }
        \toprule
        $N$ & $u$ error & $u$ EOC & $\sigma$ error & $\sigma$ EOC & $\Div_h(\sigma)$ error & $\Div_h(\sigma)$ EOC\\
        \midrule
          $2^0$ &  1.05$\times 10^{-1}$   &     --       &      6.53$\times 10^{-1}$   &     --         &        8.32$\times 10^{-1}$    &          --\\
        $2^1$ &  2.78$\times 10^{-2}$    &   1.91   &    3.07$\times 10^{-1}$     &   1.09     &      2.68$\times 10^{-1}$      &        1.64\\
        $2^2$ &  7.05$\times 10^{-3}$   &   1.98   &    1.49$\times 10^{-1}$     &   1.04       &      7.15$\times 10^{-2}$     &        1.91\\
        $2^3$ &  1.77$\times 10^{-3}$  &    2.00  &     7.38$\times 10^{-2}$   &    1.02    &       1.82$\times 10^{-2}$    &         1.98\\
        $2^4$ &  4.42$\times 10^{-4}$ &    2.00  &     3.66$\times 10^{-2}$    &    1.01    &       4.56$\times 10^{-3}$   &         1.99\\
        $2^5$ &  1.11$\times 10^{-4}$ &    2.00  &     1.83$\times 10^{-2}$   &    1.00    &       1.14$\times 10^{-3}$   &         2.00\\
        \bottomrule
    \end{tabular}
    }
    \caption{Errors and convergence rates with AW\textsuperscript{nc} near the incompressible limit $\nu = 0.4999999, \mu = 1$.}\label{tab:last_AWnc_MMS}
    \vspace{-5mm}
\end{table}


\noindent Note that the observed order of convergence in the displacement is one higher than proved, but 2\textsuperscript{nd}-order convergence for the displacement is proved for the 3D analogue of AW\textsuperscript{nc}, the nonconforming Arnold--Awanou--Winther element~\cite{arnold2014}, by applying a duality argument in the case of full elliptic regularity $(\sigma,u)\in H^1(\Omega;\mathbb{S})\times H^2(\Omega;\mathbb{R}^d)$.

\sloppy Tables~\ref{tab:first_AWc_MMS}--\ref{tab:last_AWc_MMS} exhibit convergence behaviour of the AW\textsuperscript{c} element, applied to $u(x,y) = (-\mathrm{e}^{\sin(\frac{\pi}{2}y)}, 3\cos(\pi x))^T$, and exact stress $\sigma = \mathbb{C}\varepsilon(u)$.

\begin{table}[H]
    \vspace{-3mm}
    {\tiny
    \begin{tabular}{ p{0.5cm}|p{1.8cm}|p{1.8cm}|p{1.8cm}|p{1.5cm}|p{2cm}|p{2cm} }
        \toprule
        $N$ & $u$ error & $u$ EOC & $\sigma$ error & $\sigma$ EOC & $\Div(\sigma)$ error & $\Div(\sigma)$ EOC\\
        \midrule
        $2^0$ &  1.50$\times 10^{-1}$   &     --      &     1.16$\times 10^{1}$   &       --       &          1.47$\times 10^{0}$   &            --\\
        $2^1$ &  3.81$\times 10^{-2}$    &   1.97  &   1.60$\times 10^{-2}$    &     2.85   &        3.75$\times 10^{-1}$     &          1.97\\
        $2^2$ &  9.57$\times 10^{-3}$   &   1.99   &   2.09$\times 10^{-3}$   &     2.94   &        9.44$\times 10^{-2}$   &          1.99\\
        $2^3$ &  2.40$\times 10^{-3}$  &    2.00 &    2.65$\times 10^{-4}$ &      2.98  &         2.36$\times 10^{-2}$  &           2.00\\
        $2^4$ &  5.99$\times 10^{-4}$ &    2.00 &    3.34$\times 10^{-5}$ &      2.99  &         5.91$\times 10^{-3}$ &           2.00\\
        $2^5$ &  1.50$\times 10^{-4}$ &  2.00  &     4.19$\times 10^{-6}$  &     3.00   &     1.48$\times 10^{-3}$   &         2.00     \\
        \bottomrule
    \end{tabular}
    }
    \caption{Errors and convergence rates using AW\textsuperscript{c} with $\nu = 0.25, \mu = 1$.}\label{tab:first_AWc_MMS}
    \vspace{-3mm}
\end{table}


\begin{table}[H]
    \vspace{-4mm}
    {\tiny
    \begin{tabular}{ p{0.5cm}|p{1.8cm}|p{1.8cm}|p{1.8cm}|p{1.5cm}|p{2cm}|p{2cm} }
        \toprule
        $N$ & $u$ error & $u$ EOC & $\sigma$ error & $\sigma$ EOC & $\Div(\sigma)$ error & $\Div(\sigma)$ EOC\\
        \midrule
        $2^0$ &  1.49$\times 10^{-1}$   &     --       &    1.70$\times 10^{-1}$    &      --         &         1.47$\times 10^{0}$   &           --\\
        $2^1$ &  3.81$\times 10^{-2}$    &   1.97   &  1.99$\times 10^{-2}$     &    3.09      &       3.75$\times 10^{-1}$     &         1.97\\
        $2^2$ &  9.57$\times 10^{-3}$   &   1.99   &  2.31$\times 10^{-3}$    &    3.12     &       9.44$\times 10^{-2}$    &         1.99\\
        $2^3$  & 2.40$\times 10^{-3}$  &    2.00  &   2.80$\times 10^{-4}$   &     3.05    &        2.36$\times 10^{-2}$  &          2.00\\
        $2^4$  & 5.99$\times 10^{-4}$ &    2.00  &   3.46$\times 10^{-5}$  &     3.01    &        5.91$\times 10^{-3}$  &          2.00\\
        $2^5$  & 1.50$\times 10^{-4}$ &  2.00  & 4.32$\times 10^{-6}$ &  3.00  & 1.48$\times 10^{-3}$   &  2.00   \\
        \bottomrule
    \end{tabular}
    }
    \caption{AW\textsuperscript{c} near the incompressible limit $\nu = 0.4999999, \mu = 1$.}\label{tab:last_AWc_MMS}
    \vspace{-5mm}
\end{table}


\noindent We manifestly see robustness in the incompressible regime for both elements.
The comparative magnitudes of the errors show that the AW\textsuperscript{c} element provides a markedly more accurate approximation to the stress field.
Note also that convergence behaviour in the divergence of the stress is the same for both elements, since their divergences have the same degree; this is also consistent with the `extra' AW\textsuperscript{c} vertex basis functions being solenoidal as observed in Remark~\ref{rem:vertex_divfree}.

Having validated both AW elements, we consider a more complex example from the PhD thesis of Li~\cite{lithesis} which moreover includes traction conditions. The domain, pictured in Figure~\ref{fig:lidomain}, consists of the rectangle $\Omega = [0,3]\times[0,1]$ with three disks removed, occupied by a material assumed to be isotropic and homogeneous with $\nu = 0.2$ and Young's modulus $E = 10$:
\begin{figure}[H]
    \vspace{-4mm}
    \includegraphics[width=0.8\linewidth]{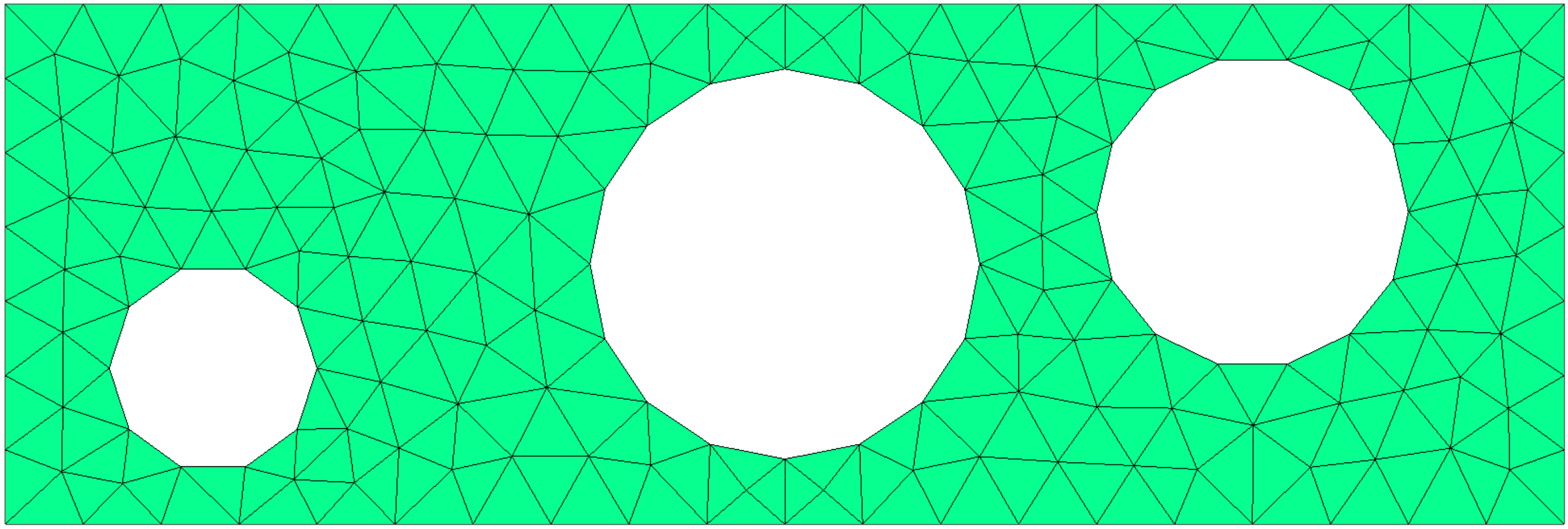}
    \vspace{-4mm}
    \caption{Domain with a coarser initial mesh than that employed in~\cite{lithesis} so that multigrid may be performed.}\label{fig:lidomain}
    \vspace{-6mm}
\end{figure}

\begin{figure}[H]
    \vspace{-3mm}
    \centering
    \includegraphics[width=0.8\linewidth]{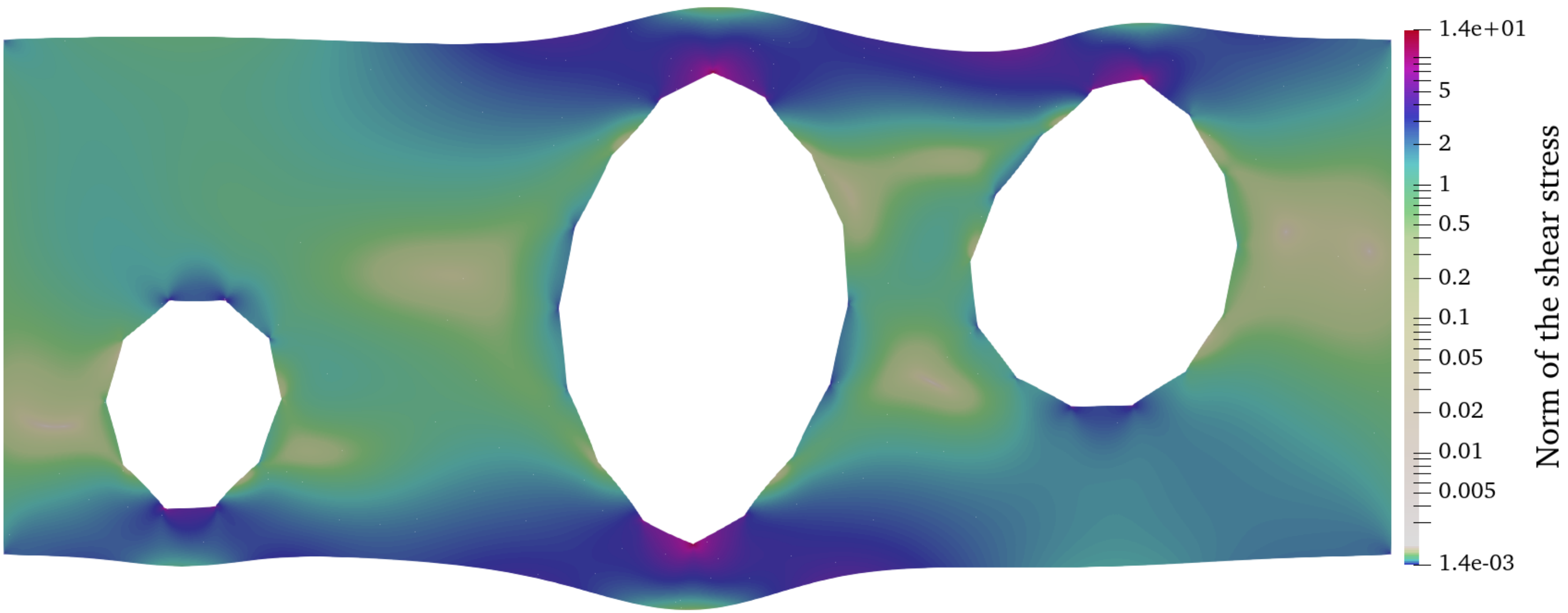}
    \vspace{-4mm}
    \caption{A traction-free condition except at both ends.}\label{fig:Li_tract}
    \vspace{-5mm}
\end{figure}

\begin{figure}[H]
    \vspace{-3mm}
    \centering
    \begin{subfigure}[t]{0.45\linewidth}
        \centering
        \includegraphics[width=1.0\linewidth]{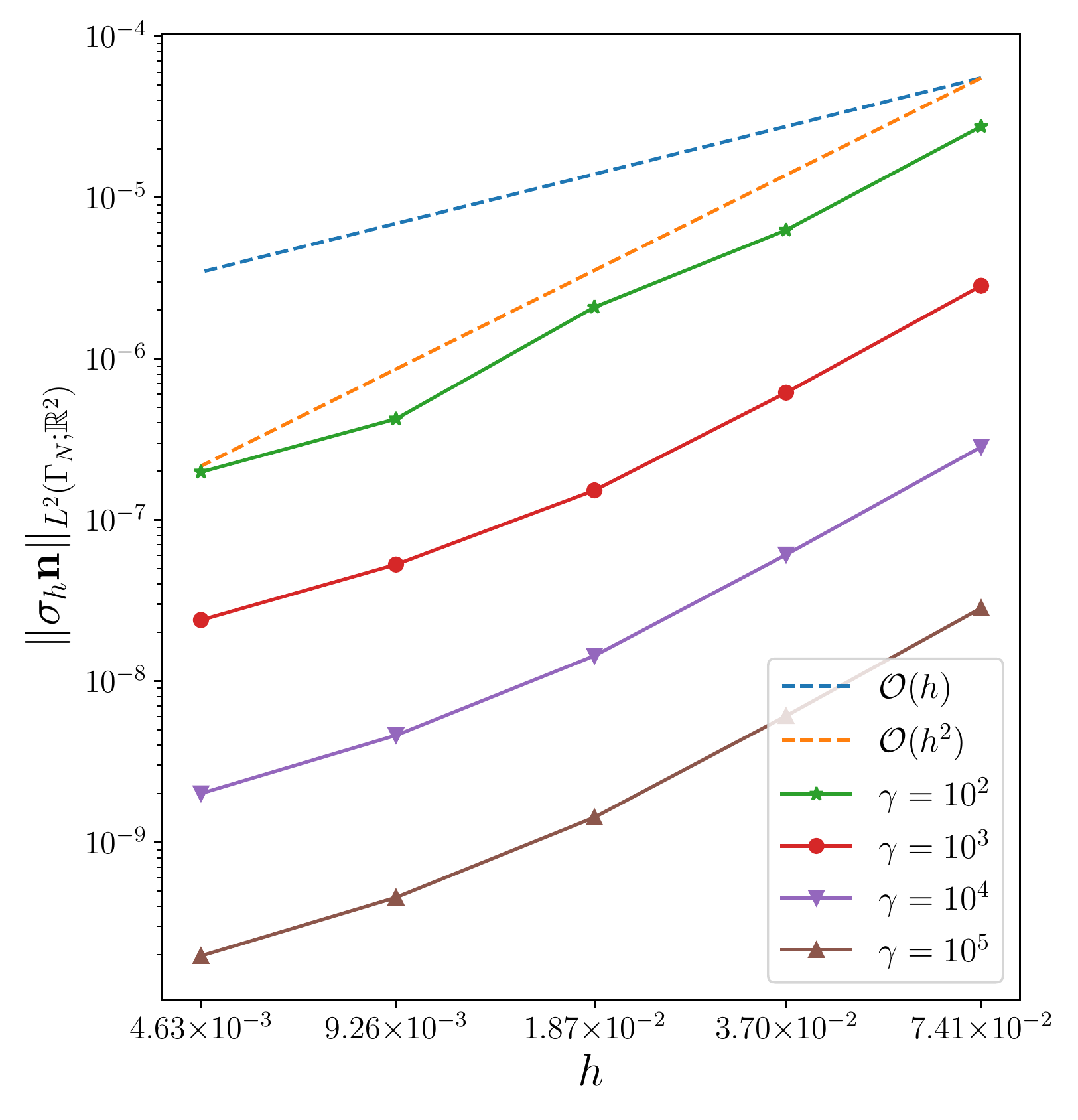}
        \vspace{-8mm}
        \caption{AW\textsuperscript{c} with $\alpha = 1$}
   \end{subfigure}
   \hspace{5mm}
   \begin{subfigure}[t]{0.45\linewidth}
        \centering
        \includegraphics[width=1.0\linewidth]{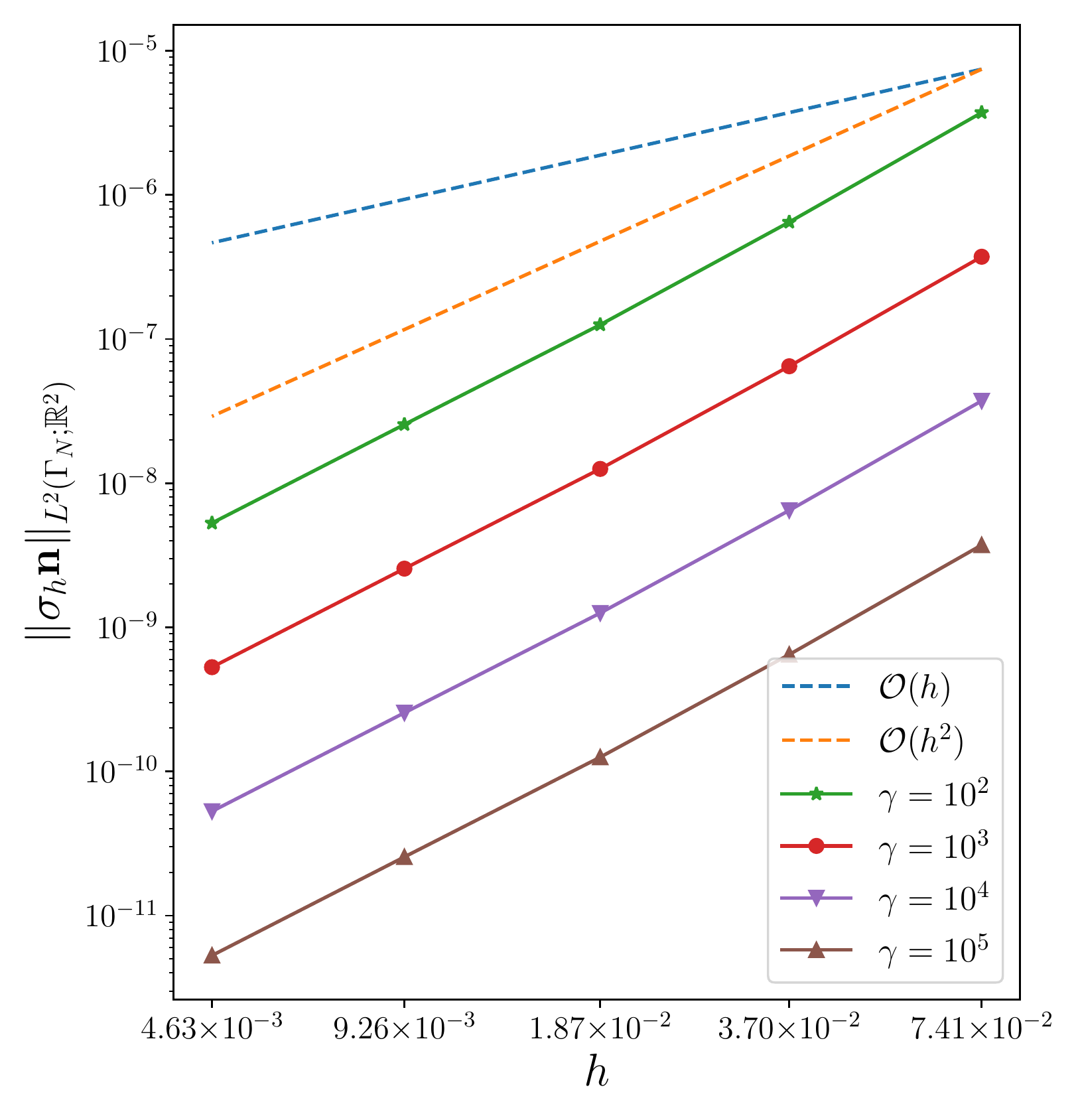}
        \vspace{-8mm}
        \caption{AW\textsuperscript{nc} with $\alpha = 0$}
    \end{subfigure}
    \vspace{-4mm}
    \caption{Convergence to the traction-free condition in $L^2(\Gamma_N;\mathbb{R}^2)$ for the conforming and nonconforming AW elements, with and without AL penalization respectively.}\label{fig:traction_conv}
    \vspace{-5mm}
\end{figure}

\noindent The prescribed displacement is fixed $(0,0)^T$ at the left-hand end and compressed $(-1,0)^T$ at the right end, which together form $\Gamma_D$, with a stress-free condition $\sigma\bn = 0$ on $\Gamma_N = \partial\Omega\setminus\Gamma_D$ given by the top, bottom, and the boundaries of the holes; there is no external force $f$. We combine the Nitsche and augmented Lagrangian penalties, seeking critical points $(\sigma_h, u_h)$ of 
\vspace{-2mm}\begin{equation}
    \mathcal{H}_{h,\gamma,\alpha}(\sigma_h, u_h) \coloneqq \mathcal{H}_{h,\gamma}(\sigma_h, u_h) + \frac{\alpha}{2}\int_\Omega \|\Div\sigma_h \|^2\dx.
\vspace{-2mm}\end{equation}
Due to the Nitsche boundary terms, which at present cannot be treated with \texttt{PCPATCH}, the application of the vertex star iteration necessitated the use of \texttt{\href{https://www.mcs.anl.gov/petsc/petsc-current/docs/manualpages/PC/PCASM.html}{\underline{PCASM}}}. 
This was carried out 
with residual $\ell^2$-norm tolerance $10^{-9}$ for the outermost Krylov iterations.
Table~\ref{tab:Li_mg} exhibits the behaviour of FGMRES for AW\textsuperscript{c}, with multigrid for the augmented stress block employing vertex star relaxation with Chebyshev smoother on each multigrid level, 
with fixed Nitsche and AL parameters $\gamma = 100, \alpha = 1$.
We verify in Figure~\ref{fig:Li_tract}, which is colored by the size of the shear stress $\|\dev\sigma_h\|$, that the AW\textsuperscript{c} solution on the finest mesh, with 
almost five million degrees of freedom, 
is free from numerical artifacts.
Figure~\ref{fig:traction_conv} shows convergence of the traction residual to zero in $L^2(\Gamma_N;\mathbb{R}^2)$ as $h\rightarrow 0$ using AW\textsuperscript{c} with $\alpha = 1$, and AW\textsuperscript{nc} (using LU factorization) with $\alpha = 0$ respectively, for various values of the Nitsche parameter $\gamma$.
Although we have proved convergence of the Nitsche penalty for any $\gamma\geq 1$, in practice we find the solver most effective at $\gamma = \mathcal{O}(100)$, and that increasing $\gamma$ further will enforce the traction condition more strongly (if desired), at the cost of a more ill-conditioned system.

\vspace{-2mm}\begin{remark}
    The observed order of convergence for the traction residual in Figure~\ref{fig:traction_conv} is higher than that predicted by Proposition~\ref{prop:nitsche_error}, which itself makes an artificial regularity assumption on the stress field, and which we therefore conjecture could be improved, for example via duality arguments.
\vspace{-2mm}\end{remark}

\begin{table}[H]
    \vspace{-3mm}
    \centering
    {\tiny
    \begin{tabular}{p{0.5cm}p{1.2cm}||c}
            \toprule
            $N$ & DOFs & Krylov iterations\\
            \midrule
            $2^1$ & 2.01$\times 10^4$ & 12\\
            $2^2$ & 7.90$\times 10^4$ & 12\\
            $2^3$ & 3.13$\times 10^5$ & 12\\
            $2^4$ & 1.25$\times 10^6$ & 12\\
            $2^5$ & 4.98$\times 10^6$ & 11\\
            \bottomrule
        \end{tabular}
        }
        \caption{Outer Krylov iteration counts using the vertex star relaxation for the AW\textsuperscript{c} element.}\label{tab:Li_mg}
        \vspace{-5mm}
\end{table}


\vspace{-3mm}\subsection{Code availability}
\vspace{-2mm}

The Mardal--Tai--Winther and conforming and nonconforming Arnold--Winther elements, with the modifications suggested by subsection~\ref{sec:scale-inv}, were incorporated into the main branch of Firedrake as part of this work.
For reproducibility, the exact software versions used to generate the numerical results in this paper are archived in~\url{https://zenodo.org/record/5596313}~\cite{zenodo/Firedrake-20211025.0}; the code, and scripts for the associated plots, are available at~\url{https://bitbucket.org/FAznaran/piola-mapped}. 

\vspace{-3mm}\section{Conclusion}
\vspace{-2mm}

We have generalized Piola transformation theory to incorporate non-standard $\hdiv$ elements, 
developing techniques for representative finite elements discretizing $\hdiv$ and $H(\Div;\mathbb{S})$ in two dimensions: the Mardal--Tai--Winther element for Stokes--Darcy flow, and two exotic, symmetry-enforcing elements for elastic stress due to Arnold and Winther. 
Numerical experiments verify the accuracy of implementations which are newly enabled by our approach. 
All elements have been implemented in the publicly available Firedrake library, within which
we have demonstrated the composability of our implementations with existing patch-based smoothers.
Of independent interest is the application of Nitsche's method to dual mixed problems; 
further investigation is merited by 
the interaction between the Nitsche and augmented Lagrangian penalties, whose efficacy when combined together we have observed numerically when applied to Arnold--Winther elements for linear elasticity.

We emphasize that our theory aims to demonstrate that unusual or non-standard elements, with desirable features but perhaps complicated transformation properties, can be used in an inexpensive, composable, and automated way, rather than to advocate for the use of MTW, AW\textsuperscript{c}, or AW\textsuperscript{nc} specifically. 
In particular, the 3D analogue of the AW spaces, namely the conforming and nonconforming Arnold--Awanou--Winther elements~\cite{arnold2008,arnold2014}, are of local dimension 42 and 162 respectively, and are presented in the understanding that ``[t]he complexity of the elements may very well limit their practical significance''~\citep[p.~1231]{arnold2008}. A cheaper alternative is provided by the conforming Hu--Zhang element, which in 3D is of dimension only 48, while many nonconforming efforts in 3D are rectangular~\cite{hu2015,yi2005}; 
in keeping with Remark~\ref{rem:H-Z}, we conjecture that the mapping properties of these 3D elements would be analogous to our analysis in Section~\ref{sec:transform}.

\vspace{-3mm}\section*{Acknowledgements}
\vspace{-2mm}
The authors are grateful to Endre S\"uli for comments on Section~\ref{sec:discretize_H-R} and to Kaibo Hu for comments on subsection~\ref{sec:pullbacks}.
\vspace{-3mm}
\bibliographystyle{plain}
\bibliography{paper}

\end{document}